\documentclass{amsart}
\usepackage{amsmath,amsthm}
\usepackage{amsfonts,amssymb}

\newtheorem{thm}{Theorem}[section]
\newtheorem{cor}[thm]{Corollary}
\newtheorem{lem}[thm]{Lemma}
\newtheorem{prop}[thm]{Proposition}

\newtheorem{defn}[thm]{Definition}
\newtheorem{re}[thm]{Remark}

\theoremstyle{remark}

\def\CL{{\mathcal L}}

\def\CS{{\mathcal S}}

\def\CV{{\mathcal V}}

\def\RR{{\mathbb R}}

\def\ZZ{{\mathbb Z}}
        
        \def\proj{\operatorname{Proj}}
        
        \def\supp{\operatorname{supp}}

\newcommand{\wh}{\widehat}

\def\cO{{\mathcal O}}
\def\cR{{\mathcal R}}
\def\cT{{\mathcal T}}
\def\cX{{\mathcal X}}

\newcommand{\eps}{{\varepsilon}}

\def\R{{R}}

\def\xxi{{\xi}}

\def\a{{\alpha}}
\def\b{{\beta}}
\def\ha{{\wh a}}

\def\W{{\mathcal W}}

\begin{document}

\title{Localized polynomial frames on the ball}

\author{Pencho Petrushev and Yuan Xu}
\address{Department of Mathematics\\University of South Carolina\\
Columbia, SC 29208.}
\email{pencho@math.sc.edu}
\address{Department of Mathematics\\ University of Oregon\\
    Eugene, Oregon 97403-1222.}\email{yuan@math.uoregon.edu}

\date{February 6, 2006}
\thanks{The first author has been supported by a grant from
the National Geospatial-Intelligence Agency.}

\begin{abstract}
Almost exponentially localized polynomial kernels are constructed
on the unit ball $B^d$ in $\RR^d$ with weights 
$W_\mu(x)= (1-|x|^2)^{\mu-1/2}$, $\mu \ge 0$, by smoothing out
the coefficients of the corresponding orthogonal projectors. These
kernels are utilized to the design of cubature formulae on $B^d$
with respect to $W_\mu(x)$ and to the construction of polynomial
tight frames in $L^2(B^d, W_\mu)$ (called needlets) whose elements
have nearly exponential localization.
\end{abstract}

\maketitle

\section{Introduction}
\setcounter{equation}{0}

The construction of bases and frames on various domains,
in particular on $\RR^d$ and on the $d$-dimensional cube, sphere, and ball,
is important from many prospectives and has numerous applications.
The example of Meyer's wavelets \cite{Meyer} and the $\varphi$-transform of
Frazier and Jawerth \cite{F-J-W} clearly shows the advantage of
using localized bases or frames for decomposition of function and distribution
spaces on $\RR^d$ in contrast to other means such as atomic decompositions
or Fourier series (in the periodic case).
Three of their features,
(i) infinite smoothness,
(ii) almost exponential space localization, and
(iii) infinitely vanishing moments,
make them a~universal tool for decomposing most of the classical
spaces on $\RR^d$, including Besov and Triebel-Lizorkin spaces.
The key to this is that the coefficients in the wavelet or
$\varphi$-transform expansions precisely capture the information
in the norms defining the corresponding spaces.

Our primary goal in this article is to develop a similar tool
for decomposition of weighted spaces of functions or distributions
on the unit ball $B^d$ in $\RR^d$ ($d>1$) with weights
$$
W_\mu(x):= (1-|x|^2)^{\mu-1/2}, \quad \mu \ge 0,
$$
where $|x|$ is the Euclidean norm of $x\in \RR^d$.
The situation here, however, is much more complicated than on $\RR^d$
(the shift invariant case) or on the torus or even on the sphere
due to several reasons:
(i) there are no dilation or translation operators on $B^d$,
(ii) the boundary of $B^d$ in combination with the weight $W_\mu(x)$
creates a great deal of inhomogeneity,
(iii) orthogonal systems such as orthogonal polynomials on $B^d$
are much less friendly than the trigonometric system, and
(iv) there are no uniformly distributed points on $B^d$ or
on the $d$-dimensional unit sphere $S^d$.

Our approach to the problem at hand will heavily rely on orthogonal polynomials
in the weighted spaces $L^2(B^d, W_\mu)$.
The standard Hilbert space theory gives the orthogonal decomposition
\begin{equation}\label{L2:decomp}
L^2(B^d, W_\mu) = \sum_{n=0}^\infty \bigoplus \CV_n^d, \qquad
     \CV_n^d \subset \Pi_n^d,
\end{equation}
where $\CV_n^d$ is the subspace of all polynomials of degree $n$ which
are orthogonal to lower degree polynomials in $L^2(B^d, W_\mu)$.
Note that $\dim \CV_n^d = \binom{n+d-1}{n}\sim n^{d-1}$, so $\CV_n^d$ is
a large subspace of $L^2$.
The orthogonal projector
$\proj_n: L^2(B^d, W_\mu) \mapsto\CV_n^d$
can be written as
$$
(\proj_n f)(x) = b_d^\mu \int_{B^d} f(y)P_n(W_\mu; x,y) W_\mu(y) dy.
$$
Here $P_n(W_\mu;x,y)$ is its kernel
and $b_d^\mu$ is the normalization
constant of $W_\mu$, namely, $(b_d^\mu)^{-1} := \int_{B^d} W_\mu(x) dx$.
It is crucial for our further development that the kernels $P_n(W_\mu;x,y)$
have an explicit representation \cite{X99} in terms of Gegenbauer polynomials
(see (\ref{eq:Kn})-(\ref{rep.Pn0}) below).
Now, (\ref{L2:decomp}) can be rewritten in the form
$$
f=\sum_{n=0}^\infty\proj_n f,
\quad f\in L^2(B^d, W_\mu).
$$
Denote by $S_n$ the orthogonal projector of $L^2(B^d, W_\mu)$ onto
$\sum_{\nu=0}^n \bigoplus \CV_\nu^d$, i.e.
$S_nf := \sum_{\nu=0}^n \proj_n f$.
Evidently,
\begin{equation}\label{def.Sn}
(S_n f)(x)= b_d^\mu \int_{B^d} f(y)K_n(W_\mu; x,y) W_\mu(y) dy
\end{equation}
with kernel
\begin{equation}\label{def.Kn}
K_n(W_\mu; x,y):= \sum_{\nu=0}^n P_\nu(W_\mu;x,y).
\end{equation}

Consider now the kernel
\begin{equation}\label{def:LnI}
L_n^\mu(x,y)=\sum_{j=0}^\infty \wh a\Big(\frac{j}{n}\Big)P_j(W_\mu;x,y),
\end{equation}
obtained by smoothing out the coefficients in the definition of the kernel
$K_n(W_\mu; x,y)$ in (\ref{def.Kn}) by sampling a smooth function $\ha$.
One of our main results in this article essentially states that if
$\ha \in C^\infty[0, \infty)$ is compactly supported, then
$L_n^\mu(x,y)$ has almost exponential (faster than any polynomial)
rate of decay away from the main diagonal $y=x$ in $B^d\times B^d$.
To state this result more precisely, let us introduce the distance
(see~(\ref{def.dist}))
\begin{equation} \label{eq:distant}
d(x,y):= \arccos \left
\{ \langle x,y\rangle + \sqrt{1-|x|^2}\sqrt{1-|y|^2} \right \}
\quad \mbox{on $B^d$}
\end{equation}
and set
$$
\W_\mu(n; x) := \left(\sqrt{1-|x|^2} + n^{-1}\right)^{2\mu}, \quad x \in B^d.
$$
Then (see \S\ref{Loc.pol.ball})
for any $k >0$ there exists a constant
$c_k>0$ depending only on $k$, $d$, $\mu$, and $\wh a$ such that
\begin{equation} \label{est.mainI}
|L_n^\mu(x,y)| \le c_k \frac{n^d }{\sqrt{\W_\mu(n;x)} \sqrt{\W_\mu(n;y)}
     (1 + n\,d(x,y))^k}.
\end{equation}

The localized kernels $L_n^\mu$ provide a powerful tool for constructing
cubature formulae on $B^d$ with weights $W_\mu(x)$, $\mu \ge 0$, that are exact
for all polynomials of degree $n$, i.e. in $\Pi_n^d$,
and have positive coefficients of the right size. It is an important feature of
our cubature formulae (see \S\ref{Cubature}) that for all $\mu \ge 0$ the knots
are obtained by projecting onto $B^d$ sets of ``almost equally" distributed
points on the upper hemisphere $S^d_+$ in $\RR^{d+1}$; the knots are in fact
almost equally distributed on $B^d$ with respect to the distance $d(\cdot,\cdot)$
defined in \eqref{eq:distant}.
Currently very few families of
cubature formulae with positive weights are known on $B^d$, among them is
the family of the product type formulae \cite{St,Myso}. However, the knots in
these formulae are not almost equally distributed.

Most importantly,
the kernels $L_n^\mu$ enable us to construct localized polynomial
frames in $L^2(B^d, W_\mu)$ which is our primary goal in this article.
Our construction is based on a semi-discrete Calder\'{o}n type
decomposition combined with our cubature formulae on the ball from
\S\ref{Cubature}.
If we denote by $\Psi=\{\psi_\xi\}_{\xi\in \cX}$ our frame on $B^d$, where
$\cX=\cup_{j=0}^\infty\cX_j$ is an index set consisting of the localization
points (poles) of the frame elements, then we have the following
representation of each $f\in L^2(B^d, W_\mu)$:
$$
f=\sum_{\xi\in \cX}\langle f, \psi_\xi\rangle\psi_\xi
\quad\mbox{and}\quad
\|f\|_{L^2(B^d, W_\mu)}
=\Big(\sum_{\xi\in \cX}|\langle f, \psi_\xi\rangle|^2\Big)^{1/2}.
$$
The above clearly indicates that $\Psi$ is a tight frame for
$L^2(B^d, W_\mu)$.
The most important feature of the frame elements $\psi_\xi$ is their
almost exponential localization:
For $\xi \in \cX_j$ (the $j$th level in $\cX$)
\begin{equation} \label{est.needlI}
|\psi_\xi(x)| \le c_k \frac{2^{jd/2}}{\sqrt{\W_\mu(2^j;x)}
(1 + 2^jd(x,y))^k}, \quad \forall k>0.
\end{equation}
Here the presence of the factor $\sqrt{\W_\mu(2^j;x)}$ is critical;
it reflects the expected influence of the boundary of $B^d$ and the
weight $W_\mu(x)$ on the localization of the frame elements.
Notice that the distance
$d(\cdot,\cdot)$ is also affected by the boundary of $B^d$.
This localization of the $\psi_\xi$'s is the reason for calling them
{\bf needlets}.
The superb localization of the needlets along with their semi-orthogonal
structure and increasing (with the levels) number of vanishing moments
enables one to utilize them for decomposition of spaces of functions
or distributions on $B^d$ other than $L^2(B^d, W_\mu)$
such as $L^p(B^d, W_\mu)$ ($1<p<\infty$) and the more general
weighted Besov and Triebel-Lizorkin spaces on $B^d$.
We will report on results of this nature in a~follow-up paper.

These ideas were first used in \cite{NPW} for the construction of
frames on the unit sphere $S^d$ in $\RR^{d+1}$.
In \cite{NPW2} the spherical frames were utilized for decomposition of
Besov and Triebel-Lizorkin spaces on the sphere.
Further, this scheme has been applied in \cite{PX} for the development
of frames on $[-1, 1]$ with Jacobi weights and
then used in \cite{KPX} for decomposition of weighted Besov and Triebel-Lizorkin
spaces on the interval.

This article is organized as follows.
In \S\ref{Principles} we outline the general principles
which guide us in constructing localized kernels and
frames on domains other than $\RR^d$.
In \S\ref{Loc.pol.int} we present some results on localized polynomial
kernels on $[-1, 1]$ with Jacobi weights.
In \S\ref{Loc.pol.ball} we prove our main results on
localized polynomial kernels on $B^d$ with weights $W_\mu(x)$, $\mu \ge 0$.
In \S\ref{Cubature} we construct cubature formulae on
$B^d$ with weights $W_\mu(x)$.
In \S\ref{Frame} we construct our needlet system and give some of its
properties.
Section \ref{appendix} is an appendix, where we give the proofs of some
results from the previous sections.

Throughout this paper positive constants are denoted by $c, c_1,\dots$
and they may vary at every occurrence.
As usual the constants may depend on some parameters, which are indicated
explicitly in some important cases.
The notation $A\sim B$ means $c_1A\le B\le c_2 A$.

\section{General principles for constructing localized kernels
and frames}\label{Principles}
\setcounter{equation}{0}

Let $(E, \mu)$ is a measure space with $E$ a metric space
and suppose that there is
an orthogonal decomposition of $L^2(E, \mu)$:
\begin{equation}\label{orth-decomp}
L^2(E, \mu) = \sum_{n=0}^\infty \bigoplus \CV_n,
\end{equation}
where $\CV_n$ is a subspace of dimension $\dim \CV_n\sim n^\gamma$, $\gamma>0$.
Let $P_n$ be the kernel of the orthogonal projector
$\proj_n: L^2(E, \mu)\to\CV_n$, i.e.
$$
(\proj_n f)(x) =\int_E P_n(x,y)f(y)d\mu, \quad f\in L^2(E, \mu).
$$
Notice that $P_n$ can be written in the form
$
P_n(x, y)=\sum_{j=1}^{\dim \CV_n}p_j(x)\overline{p_j(y)},
$
where $\{p_j\}$ is an orthonormal basis for $\CV_n$.
Then $K_n:=\sum_{j=0}^n P_\nu$ is the kernel of the orthogonal projector
onto $\sum_{\nu=0}^n \bigoplus \CV_\nu$.
In most cases of interest the kernel $K_n(x,y)$ has poor localization,
examples include the trigonometric system, orthogonal polynomials in
one or several variables on various domains.

\smallskip

\noindent
{\bf Localization principle.}
It seems to us that there is a general localization principle, which says that
for all ``natural" orthogonal systems, if the coefficients of the kernel
$K_n$ are smoothed out as in (\ref{def:LnI}) by sampling a $C^\infty$
function, then the resulting kernel has ``excellent" localization
around the main diagonal $y=x$ in $E\times E$.
To be more specific, suppose $\ha \in C^\infty(\RR)$, $\ha$ is even, and
$\ha$ is compactly supported or
$\ha \in \CS$ (the Schwartz class of rapidly decreasing $C^\infty$
functions on $\RR$). Define
\begin{equation} \label{def:Ln-gen}
L_n(x,y) := \sum_{j=0}^\infty \wh a\Big(\frac{j}{n}\Big) P_j(x,y).
\end{equation}
Then for all "natural" orthogonal systems, the kernel
$L_n(x,y)$ decays away from the main diagonal $y=x$ at nearly exponential
(faster than any polynomial) rate with respect to the distance in $E$.


In the case of the trigonometric system this principle is well-known
and widely used.
It is a fundamental fact in Harmonic Analysis that the Fourier
transform of every function $f$ in
the Schwartz space $\CS=\CS(\RR^d)$ belongs to the same space.
Consequently, any trigonometric polynomial
$ L_n(t) := \sum_{\nu=-n}^n a_\nu e^{i\nu t}$
with coefficients $\{a_\nu\}$ obtained by sampling
a compactly supported $C^\infty$ function has
faster than any polynomial rate of decay away from zero.
To make this more precise, let
$$
L_n(t) :=\sum_{\nu\in\ZZ} \ha\Big(\frac{\nu}{n}\Big)e^{i\nu t},
$$
where $\ha$ is compactly supported and $\ha \in C^\infty(\RR)$.
Then $L_n$ is a trigonometric polynomial of degree $cn$
and one easily shows the following localization of $L_n$:
For any $k > 0$ and $r\ge 0$ there exists a constant $c_k>0$
depending only on $k$, $r$, and $\ha$ such that
\begin{equation}
\label{L.bound} |L_n^{(r)}(t)|\le c_k\frac{n^{r+1}}{(1+n|t|)^k},
\quad t\in [-\pi,\pi].
\end{equation}
This estimate will serve as a prototype for our further localization
results.

For Gegenbauer polynomials and spherical harmonics
the localization principle is established and used in \cite{NPW}
and also follows by the general result in \cite{H} on the spectral
properties of elliptic operators.
For Jacobi polynomials it is proved in \cite{BD} and \cite{PX}
(see Theorem~\ref{thm:loc.int} below).
For Hermite and Laguerre polynomials the localization principle is
established in \cite{Dz}.
We will establish it here for multivariate orthogonal polynomials
in $L^2(B^d, W_\mu)$ (see Theorem~\ref{thm:loc.ball}).
We believe that the localization principle is valid in more general
settings as well.

For our purposes we restrict our attention to
``smoothing functions" $\wh a$ satisfying:


\begin{defn} \label{defn:admissible}
A function $\wh a$ is said to be admissible if
$\wh a \in C^\infty[0, \infty)$, $\wh a(t) \ge 0$, and
$\wh a$ satisfies one of the following two conditions:
\begin{enumerate}
\item[(a)] $\supp \wh a \subset [0,2]$, $\wh a(t) =1$ on $[0, 1]$, and
$0 \le \wh a (t) \le 1$ on $[1,2]$; or
\item[(b)] $\supp \wh a \subset [1/2,2]$.
\end{enumerate}
\end{defn}

There are two important applications of the localized kernels $L_n(x,y)$
from~(\ref{def:Ln-gen}):

(i) If $\ha$ is admissible of type (a), then the operator
$$
(\CL_n f)(x):=\int_E L_n(x,y)f(y)d\mu(y)
$$
apparently satisfies:
$\CL_n f=f$ for all $f \in \sum_{\nu=0}^n \bigoplus \CV_\nu$
and $\CL_n f \in \sum_{\nu=0}^{2n} \bigoplus \CV_{\nu}$.
These along with the superb localization of $L_n$ (to be established)
makes $\CL_n$ a~useful tool.
We will see this operator at work in the construction of cubature
formulae on the ball in \S \ref{Cubature}.

(ii) Kernels $L_n(x,y)$ with $\ha$ admissible of type (b)
are a valuable tool for constructing localized frames.
Let, in addition, $\wh a$ satisfy the conditions:
$\wh a(t) \ge 0$ and
\begin{equation}\label{a3G}
\wh a^2(t) + \wh a^2(2t) =1,
\quad t \in [1/2, 1].
\end{equation}
Then
\begin{equation}\label{a4G}
\sum_{\nu=0}^\infty \wh a^2(2^{-\nu}t)= 1,
\quad t \in [1, \infty).
\end{equation}
It is easy to construct such functions (see e.g. \cite{NPW}).
Define
\begin{equation}\label{def.LambdajG}
L_0(x,y) := P_0(x,y) \quad \mbox{and} \quad
L_j(x, y) := \sum_{\nu=0}^\infty
\wh a \Big(\frac{\nu}{2^{j-1}}\Big)P_\nu(x,y),
\quad  j=1, 2, \dots,
\end{equation}
and denote briefly
$$
(L_j*f)(x):= \int_E L_j(x, y)f(y)d\mu(y).
$$
One easily obtains the following semi-discrete
Calder\'{o}n type decomposition (see e.g. \cite{PX})
\begin{equation}\label{L.reprG}
f=\sum_{j=0}^\infty L_j*L_j*f
\quad \mbox{for $ f\in L_2(E, \mu)$.}
\end{equation}
To get a completely discretized decomposition of $L_2(E, \mu)$
one can use quadrature (cubature) formulae, if available.
Assume that there is a quadrature formula
\begin{equation}\label{quadrature}
\int_E f d\mu \sim \sum_{\xxi \in \cX_j} \lambda_\xxi f(\xxi)
\end{equation}
with $\cX_j \subset E$ and $\lambda_\xxi >0$,
which is exact for all functions $f$
of the form $f=gh$ with
$g, h \in \sum_{\nu=0}^{2^{2j}} \bigoplus \CV_{\nu}$.

After these preparations we now define the frame elements by
\begin{equation}\label{def.frameG}
\psi_\xi(x):= \sqrt{\lambda_\xi}\cdot L_j(x, \xi)
\quad \mbox{for}\quad
\xi\in \cX_j, ~ j=0, 1, \dots.
\end{equation}
The $\psi$'s inherit the localization of the kernels $L_j$,
which is almost exponential in all cases of interest.
This is the reason for calling them needlets.

We write $\cX := \cup_{j = 0}^\infty \cX_j$, where any two points
$\xxi, \omega \in \cX$ (from levels $\cX_j\ne \cX_k$)
are considered to be different elements of $\cX$ even if they coincide.
We use $\cX$ as an index set in the definition of the needlet system
$$
\Psi:=\{\psi_\xi\}_{\xi\in\cX}.
$$

The next statement shows that
$\Psi$ is a tight frame in $L^2(E, \mu)$.


\begin{prop}\label{prop:frameG}
If $f\in L^2(E, \mu)$, then
\begin{equation}\label{frame1G}
f =\sum_{j=0}^\infty\sum_{\xi\in\cX_j}
  \langle f, \psi_{\xi} \rangle \psi_{\xi}
= \sum_{\xi\in \cX} \langle f, \psi_\xi\rangle \psi_{\xi}
\quad\mbox{in $L^2(E, \mu)$}
\end{equation}
and
\begin{equation}\label{frame2G}
\|f\|_{L^2(E, \mu)}
=\Big(\sum_{\xi\in \cX} |\langle f, \psi_\xi\rangle|^2\Big)^{1/2}.
\end{equation}
\end{prop}
The proof of this proposition is a mere repetition of the proof
of Theorem 4 in \cite{PX} (see also \cite{NPW}) and will be omitted.

\section{Localized polynomial kernels on $[-1, 1]$}
\label{Loc.pol.int}
\setcounter{equation}{0}

The Jacobi polynomials $\{P_n^{(\a,\b)} \}_{n=0}^\infty$
constitute an orthogonal basis for the weighted space
$L^2([-1, 1], w_{\a,\b})$ with $w_{\a,\b}(t):=(1-t)^\a(1+t)^\b$,
$\a, \b>-1$.
We let $c_{\a,\b}$ denote the normalization constant of
$w_{\a,\b}$, i.e.
$c_{\a,\b}^{-1} := \int_{-1}^{1} w_{\a,\b}(t)dt$.
It is well known that~\cite{Sz}
$$
  c_{\a,\b} \int_{-1}^1 P_n^{(\a,\b)}(t)
    P_m^{(\a,\b)}(t)w_{\a,\b}(t)dt  = \delta_{n, m} h_n^{(\a,\b)},
$$
where
$$
h_n^{(\a,\b)} = \frac{\Gamma(\a+\b+2)}{\Gamma(\a+1)\Gamma(\b+1)}
\frac{\Gamma(n+\a+1)\Gamma(n+\b+1)}{(2n+\a+\b+1)\Gamma(n+1)\Gamma(n+\a+\b+1)}.
$$
For $f \in L^2([-1, 1], w_{\a,\b})$ the Fourier expansion of $f$ in
Jacobi polynomials is
$$
  f(t) = \sum_{n=0}^\infty d_n(f) (h_n^{(\a,\b)})^{-1}
    P_n^{(\a,\b)}(t),
  \quad d_n(f) = c_{\a,\b}\int_{-1}^1 f(t) P_n^{(\a,\b)}(t)
            w_{\a,\b}(t)dt.
$$
The $n$th partial sum of this expansion can be written as
$$
  (S_n f)(x) = \sum_{j=0}^n d_j(f) (h_j^{(\a,\b)})^{-1}
     P_j^{(\a,\b)}(x)
   = c_{\a,\b} \int_{-1}^1 f(t) K_n^{(\a,\b)}(x,t) w_{\a,\b}(t)dt,
$$
where the kernel is given by
\begin{equation}\label{eq:1.1}
K_n^{(\a,\b)}(x,t) = \sum_{j=0}^n
  \Big(h_j^{(\a,\b)}\Big)^{-1} P_j^{(\a,\b)}(x) P_j^{(\a,\b)}(y).
\end{equation}


The grand question here is: What is the localization
around the main diagonal $y=x$ in $[-1, 1]^2$
of a polynomial kernel of the form
\begin{equation}\label{def.L}
L_n^{\a,\b}(x,y)=\sum_{j=0}^\infty \wh a \Big(\frac{j}{n}\Big)
     \Big(h_j^{(\a,\b)}\Big)^{-1} P_j^{(\a,\b)}(x) P_j^{(\a,\b)}(y),
\end{equation}
where $\wh a \in C^\infty$?

To address this question, denote
$$
w_{\a,\b}(n; x) := (1-x + n^{-2})^{\a+1/2} (1+x +
    n^{-2})^{\b+1/2}. 
$$


\begin{thm}\label{thm:loc.int}\cite{PX}
Let $\alpha, \beta > -1/2$ and let $\wh a$
be admissible according to Definition~\ref{defn:admissible}.
Then for every $k >0$ there is a constant $c_k > 0$ depending only on
$k$, $\a$, $\b$, and $\wh a$ such that
for $0 \le \theta,\phi \le \pi$

\begin{equation}\label{Lbound1}
|L_n^{\a,\b} (\cos\theta, \cos \phi)|
\le c_k \frac{n}{\sqrt{w_{\a,\b}(n; \cos \theta)}
\sqrt{w_{\a,\b}(n; \cos \phi)}(1+n|\theta - \phi|)^{k}}.
\end{equation}
Here the dependence of $c_k$ on $\wh a$ is of the form
$c_k=c(\a,\b,k)\max_{0\le\nu\le k}\|\wh a^{(\nu)}\|_{L^\infty}$.
\end{thm}

For the proof of this theorem it is important to establish
estimate (\ref{Lbound1}) first in the particular case when $\phi=0$
(the localization of $L_n^{\a,\b}(x,1)$).
Set
\begin{equation}\label{def.Ln1}
L_n^{\a,\b}(x) := L_n^{\a,\b}(x,1)
=\sum_{j=0}^\infty \wh a\Big(\frac{j}{n}\Big)
  \Big(h_j^{(\a,\b)}\Big)^{-1} P_j^{(\a,\b)}(1)P_j^{(\a,\b)}(x).
\end{equation}
Since \cite[(4.1.1), p. 58]{Sz}
$$
P_n^{(\a,\b)}(1) =\binom{n+\a}{n}
=\frac{\Gamma(n+\a+1)}{\Gamma(\a+1)\Gamma(n+1)},
$$
it is easy to verify that
\begin{equation}\label{eq:L-n}
L_n^{\a,\b}(x) =  c^\diamond
\sum_{j=0}^\infty \wh a\Big(\frac{j}{n}\Big)
\frac{(2j+\a+\b+1)\Gamma(j+\a+\b+1)}{\Gamma(j+\b+1)}P_j^{(\a,\b)}(x),
\end{equation}
where
$c^\diamond:=\Gamma(\b+1)/\Gamma(\a+\b+2)$.

Now the key role is played by the following theorem,
which will also be critical for the proof of our main localization
result (Theorem~\ref{thm:loc.ball}).


\begin{thm} \label{thm:BD}\cite{BD, PX}
Let $\wh a$ be admissible and
assume that $\a \ge \b > -1/2$.
Then for every $k>0$ and $r\ge 0$
there exists a constant $c_k>0$ depending only on $k$, $r$, $\a$, $\b$,
and $\wh a$ such that
\begin{equation}\label{est.Ln}
\Big|\frac{d^r}{dx^r} L_n^{\a,\b}(\cos \theta)\Big|
\le c_k \frac{n^{2\a+2r+2}}{(1+n \theta)^{k}},
\quad 0 \le \theta \le \pi.
\end{equation}
The dependence of $c_k$ on $\wh a$ is of the form
$c_k=c(\a,\b,k,r)\max_{0\le\nu\le k}\|\wh a^{(\nu)}\|_{L^\infty}$.
\end{thm}

This theorem is proved in
\cite{BD} with $\ha$ admissible of type (a)
and in \cite{PX} with $\ha$ admissible of type (b).
The proof in \cite{PX} rests on the localization properties
of trigonometric polynomials given in (\ref{L.bound}),
while the proof in \cite{BD} is based on a property of Jacobi polynomials;
it can be carried out with $\ha$ admissible of type (b) as well.
Estimate (\ref{est.Ln}) was proved earlier in \cite{NPW} in the case
$\a=\b=\lambda-1/2$ (with $\lambda$ a half integer)
and utilized for the construction of frames on
the $n$ dimensional sphere.
For the reader's convenience we give the proof of Theorem~\ref{thm:BD}
(following the idea from \cite{BD}) in the appendix.

Theorem~\ref{thm:loc.int} is established in \cite{PX}.
Its proof rests on Theorem~\ref{thm:BD}.

\section{Localized polynomial kernels on the unit ball}
\label{Loc.pol.ball}
\setcounter{equation}{0}

It is known (see \cite{X99}) that the orthogonal projector
$\proj_n: L^2(B^d, W_\mu) \mapsto\CV_n^d$
can be written as
$$
(\proj_n f)(x) = b_d^\mu \int_{B^d} f(y)P_n(W_\mu; x,y) W_\mu(y) dy,
$$
where if $\mu>0$ the kernel $P_n(W_\mu;x,y)$ has the following explicit
representation:
\begin{align} \label{eq:Kn}
P_n(W_\mu; x,y)& = b_1^{\mu-\frac{1}{2}}\frac{\lambda+n}{\lambda} \\
&\times  \int_{-1}^1 C_n^\lambda \left(\langle x,y\rangle +
u\sqrt{1-|x|^2} \sqrt{1-|y|^2}\right) (1-u^2)^{\mu-1}du,  \notag
\end{align}
where $\langle x, y \rangle$ is the usual Euclidean inner product,
$C_n^\lambda$ is the $n$th degree Gegenbauer polynomial,
and
$$
   \lambda = \mu + \frac{d-1}{2}.
$$
The case $\mu = 0$ is a limit case and we have
\begin{align}\label{rep.Pn0}
 P_n(W_0; x,y) = \frac{\lambda+n}{2 \lambda} & \left[  C_n^\lambda
 \left(\langle x,y\rangle +\sqrt{1-|x|^2} \sqrt{1-|y|^2}\right) \right .\\
& \left . + C_n^\lambda \left (\langle x,y\rangle -
\sqrt{1-|x|^2} \sqrt{1-|y|^2} \right)\right].
  \notag
\end{align}

For an admissible $\wh a$ (according to Definition~\ref{defn:admissible})
we define
\begin{equation*} 
L_n^\mu(x,y)=\sum_{j=0}^\infty \wh a\Big(\frac{j}{n}\Big)P_j(W_\mu;x,y),
\quad x,y \in B^d.
\end{equation*}
The explicit representation \eqref{eq:Kn} gives
\begin{equation} \label{rep.Ln}
L_n^\mu(x,y) = b_1^{\mu-\frac{1}{2}}
  \int_{-1}^1 L_n^\lambda (\langle x,y\rangle +
    u\sqrt{1-|x|^2} \sqrt{1-|y|^2}) (1-u^2)^{\mu-1}du,
\end{equation}
where $L_n^\lambda$ is defined by
$$
 L_n^\lambda(x) = \sum_{j=0}^\infty \wh a\Big(\frac{j}{n}\Big)
          \frac{j+\lambda}{\lambda} C_j^\lambda (x).
$$
Since
$$
C_n^\lambda(x)=\frac{\Gamma(\lambda+1/2)}{\Gamma(2\lambda)}
\frac{\Gamma(n+2\lambda)}{\Gamma(n+\lambda+1/2)}
P_n^{(\lambda-1/2, \lambda-1/2)}(x),
$$
it readily follows from (\ref{eq:L-n}) that
$
L_n^{\lambda-1/2, \lambda-1/2}(x)=L_n^\lambda(x).
$
Then by Theorem~\ref{thm:BD} we get the following estimate:
For all $k, \lambda >0$ and $r\ge 0$ there exists a constant
$c_k>0$ depending only on $k$, $r$, $\lambda$, and $\wh a$, such that
\begin{equation}\label{est-Llam}
\Big|\frac{d^r}{dx^r}L_n^\lambda(\cos \theta)\Big|
\le c_k \frac{n^{2\lambda+2r+1}}{(1+n \theta)^k},
\quad 0 \le \theta \le \pi.
\end{equation}

\smallskip
\noindent
{\bf Distance on \boldmath $B^d$.}
In order to show that
$L_n^\mu$ is a well localized kernel and for our further development,
we need to introduce an appropriate distance in $B^d$ that
takes into account the fact that $B^d$ has a boundary.
In \cite{X98} it is shown that the orthogonal polynomials on the unit
ball and those on the unit sphere are closely related by the simple map
\begin{equation}\label{eq:lift}
x \in B^d \mapsto x':= (x, \sqrt{1-|x|^2}) \in S^d,
\end{equation}
which ``lifts" the points from $B^d$ to the upper hemisphere $S^d_+$
in $\RR^{d+1}$, that is, $S^d_+:=\{x\in S^{d}:x_{d+1}\ge 0\}$.
This relation leads us to
the following distance on $B^d$, which will play a vital role in the
following:
\begin{equation}\label{def.dist}
d(x,y):= \arccos \left
  \{ \langle x,y\rangle + \sqrt{1-|x|^2}\sqrt{1-|y|^2} \right \}.
\end{equation}
In fact this is the geodesic distance between
$x':= (x, \sqrt{1-|x|^2})$ and $y' := (y, \sqrt{1-|y|^2})$
on $S^d_+ \subset \RR^{d+1}$
and, consequently, it is a true distance on $B^d$. This distance has
been used to prove various polynomial inequalities, see the discussions
in \cite{BLW} and the references therein.

The map \eqref{eq:lift} also leads to a close relation between the spaces
$L^p(B^d, W_0)$ and $L^p(S^d, d\omega)$, where $d\omega$ is the surface
measure on $S^d$.
This allows us to derive results on $L^p(B^d, W_0)$ from
those on $L^p(S^d, d\omega)$, which are also easier to prove.
For these reasons we will prove our results only in the case $\mu > 0$.

The following lemma provides an important relation between $d(\cdot, \cdot)$
and the Euclidean norm $|\cdot|$ in $B^d$.


\begin{lem}\label{lem:norm-dist}
For $x, y\in B^d$, we have
\begin{equation} \label{norm-dist1}
\big||x|-|y|\big|\le \frac{1}{\sqrt{2}}d(x, y)
  \left(\sqrt{1-|x|^2}+\sqrt{1-|y|^2}\right)
\end{equation}
and hence
\begin{equation} \label{norm-dist2}
\Big|\sqrt{1-|x|^2}-\sqrt{1-|y|^2}\Big| \le \sqrt{2}\, d(x, y).
\end{equation}
\end{lem}

\begin{proof}
Let $0\le \alpha, \beta \le \pi/2$ be defined from $|x|=\cos \alpha$
and $|y|=\cos \beta$. Using spherical-polar coordinates $x = |x| \xi$ and
$y = |y| \zeta$, where $\xi,\zeta \in S^{d-1}$, we see that
$$
 d(x,y) = \arccos \left( \cos \alpha \cos\beta \langle \xi, \zeta \rangle +
    \sin \alpha \sin \beta \right) \ge \arccos (\cos(\alpha - \beta))
$$
which yields $d(x,y) \ge |\alpha -\beta|$.
On the other hand, since
$0 \le \alpha,\beta \le \pi/2$, we have $\cos \frac{\alpha-\beta}{2} \ge
\cos(\pi/4)=\sqrt{2}/2$ and, consequently,
$$
\sin \alpha + \sin \beta = 2 \sin \frac{\alpha+\beta}{2}
    \cos \frac{\alpha-\beta}{2} \ge \sqrt{2} \sin \frac{\alpha+\beta}{2}.
$$
Using the above we obtain
\begin{align*}
\big||x|-|y|\big| &= |\cos \alpha - \cos\beta|
 = 2\sin\frac{|\alpha-\beta|}{2}\sin\frac{\alpha+\beta}{2}\\
&\le \frac{1}{\sqrt{2}}|\alpha-\beta|(\sin\alpha+\sin\beta)
 \le \frac{1}{\sqrt{2}} d(x, y)(\sqrt{1-|x|^2}+\sqrt{1-|y|^2}).
\end{align*}
Thus (\ref{norm-dist1}) is established. Estimate (\ref{norm-dist2}) follows
immediately from (\ref{norm-dist1}).
\end{proof}

Let us define
\begin{equation}\label{def.W}
\W_\mu(n; x) := \left(\sqrt{1-|x|^2} + n^{-1}\right)^{2\mu}, \quad x \in B^d.
\end{equation}
Our next theorem shows that the kernels $L_n^\mu$
are almost exponentially localized around the main diagonal $y=x$
in $B^d\times B^d$.


\begin{thm} \label{thm:loc.ball}
Let $\wh a$ be admissible. Then for any $k >0$ there exists a constant
$c_k>0$ depending only on $k$, $d$, $\mu$, and $\wh a$ such that
\begin{equation} \label{est.main}
|L_n^\mu(x,y)| \le c_k \frac{n^d }{\sqrt{\W_\mu(n;x)} \sqrt{\W_\mu(n;y)}
     (1 + n\,d(x,y))^k},
\quad x, y\in B^d.
\end{equation}
\end{thm}


\begin{re} \label{re:Ck}
Theorem~\ref{thm:loc.ball} as well as Theorems~\ref{thm:loc.int}-\ref{thm:BD}
can obviously be modified for the case when $\ha\in C^k$
$($$k$ sufficiently large$)$ in place of $\ha\in C^\infty$.
\end{re}

We will derive Theorem~\ref{thm:loc.ball}
when $\mu>0$ from estimate (\ref{est-Llam})
and the following lemma, using representation (\ref{rep.Ln}) of $L_n^\mu$.
The proof in the case $\mu=0$ is easier
and will be omitted; it utilizes (\ref{rep.Pn0}).

Let us denote briefly
\begin{equation} \label{def.t}
t(x,y;u):=\langle x, y\rangle + u\sqrt{1-|x|^2}\sqrt{1-|y|^2}.
\end{equation}


\begin{lem} \label{lem:loc1}
Let $\gamma >-1$, $k >3\gamma+4$, and $n\ge 1$.
Then for $x,y \in B^d$
\begin{align}\label{loc.est}
&\int_{-1}^1\frac{(1-u^2)^\gamma du}{(1+n\sqrt{1-t(x,y;u)})^k}\\
&\le \frac{cn^{-2\gamma-2}}
{(\sqrt{1-|x|^2}+n^{-1})^{\gamma+1}(\sqrt{1-|y|^2}+n^{-1})^{\gamma+1}
(1+nd(x,y))^{k-3\gamma-4}},\notag
\end{align}
where $c>0$ depends only on $\gamma$, $k$, and $d$.
\end{lem}

\begin{proof}
Denote briefly
$t:=t(x,y;u)$.
Then we can write
\begin{align*}
1-t   =  1- \langle x, y\rangle - \sqrt{1-|x|^2}\sqrt{1-|y|^2}
    + (1-u)\sqrt{1-|x|^2}\sqrt{1-|y|^2},
\end{align*}
which implies
\begin{align}\label{eq:2.8}
   1- t & \ge 1 - \langle x, y\rangle - \sqrt{1-|x|^2}\sqrt{1-|y|^2} \\
     & = 1- \cos d(x,y) = 2 \sin^2 \frac{d(x,y)}{2}
   \ge \frac{2}{\pi^2} d(x,y)^2 \notag
\end{align}
and
\begin{align}\label{eq:2.9}
1 - t &\ge \frac{2}{\pi^2} d(x,y)^2+(1-u)\sqrt{1-|x|^2}\sqrt{1-|y|^2}\\
&\ge (1-u)\sqrt{1-|x|^2}\sqrt{1-|y|^2}.\notag
\end{align}
By (\ref{eq:2.8}), we have
\begin{align}\label{est-L1}
\int_{-1}^1\frac{(1-u^2)^\gamma du}{(1+n\sqrt{1-t})^k}
  \le   \frac{c}{(1+n d(x,y))^k}.
\end{align}
Inequality (\ref{loc.est}) will follow from this and the estimate:
\begin{align} \label{eq:3.4}
& \int_{-1}^1\frac{(1-u^2)^\gamma du}{(1+n\sqrt{1-t})^k} \\
& \qquad\qquad \le \frac{cn^{-2\gamma-2}}
{(\sqrt{1-|x|^2})^{\gamma+1}(\sqrt{1-|y|^2})^{\gamma+1}
(1+n d(x,y))^{k-2\gamma-3}}. \notag
\end{align}
To establish this last estimate, we split the integral over $[-1,1]$ into two
integrals: one over $[-1,0]$ and the other over $[0,1]$.
For the integral over $[-1,0]$ we write the factor $(1+n\sqrt{1-t})^k$
as the product of
$(1+n\sqrt{1-t})^{k-2\gamma-2}$ and
$(1+n\sqrt{1-t})^{2\gamma+2}$.
Then we apply inequalities \eqref{eq:2.8} and \eqref{eq:2.9}
to the first and the second terms, respectively.
This gives
\begin{align*}
\int_{-1}^0 &\le
    \frac{c}{(1+n d(x,y))^{k-2\gamma-2}}
  \int_{-1}^0 \frac{(1-u^2)^{\gamma} }
  { [n^2 \sqrt{1-|x|^2} \sqrt{1-|y|^2}(1-u)]^{\gamma+1}} du   \\
&\le \frac{ cn^{-2\gamma-2}}
{(\sqrt{1-|x|^2} )^{\gamma+1}(\sqrt{1-|y|^2} )^{\gamma+1}
          (1+n d(x,y))^{k-2\gamma-2}}.
\end{align*}
We now estimate the integral over $[0,1]$.
Denote briefly $A: = \sqrt{1-|x|^2} \sqrt{1-|y|^2}$.
Using (\ref{eq:2.9}) and applying the substitution $s=An^2(1-u)$, we get
\begin{align*}
\int_0^1 &\le
c\int_0^1 \frac{ (1-u^2)^{\gamma}}
{(1+n \sqrt{d(x,y)^2+ A(1-u)})^k} du \\
&\le \frac{c}{(An^2)^{\gamma+1}}
  \int_0^{An^2} \frac{ s^{\gamma} }
  {(1+ \sqrt{n^2 d(x,y)^2+ s})^k}ds \\
&\le  \frac{cn^{-2\gamma-2}}{A^{\gamma+1}
   (1+n d(x,y))^{k-2\gamma-3} } \int_0^\infty
    \frac{ s^{\gamma} ds }
      {(1+ \sqrt{n^2 d(x,y)^2+ s})^{2\gamma+3}}\\
&\le  \frac{cn^{-2\gamma-2}}{A^{\gamma+1}
         (1+n d(x,y))^{k-2\gamma-3}}.
\end{align*}
Putting these estimates together gives \eqref{eq:3.4}.

To complete the proof of (\ref{loc.est}) we need the following
simple inequality (see inequality (2.21) in \cite{PX}):
\begin{equation} \label{est.ab}
(a + n^{-1})(b+n^{-1}) \le 3(ab+n^{-2})(1+n|a-b|),
\quad a,b \ge 0, \: n\ge 1.
\end{equation}
Inequalities (\ref{norm-dist2}) and (\ref{est.ab}) yield
\begin{align*}
&(\sqrt{1-|x|^2}+n^{-1})(\sqrt{1-|y|^2}+n^{-1})\\
&\qquad\qquad\le 3(\sqrt{1-|x|^2}\sqrt{1-|y|^2}+n^{-2})
(1+nd(x,y)).
\end{align*}
This along with (\ref{est-L1}) and (\ref{eq:3.4}) implies
(\ref{loc.est}).
\end{proof}


\noindent
{\em Proof of Theorem \ref{thm:loc.ball}.}
For $t = \cos \theta$, $0 \le \theta \le \pi$, we have $\theta /2 \sim
\sin \theta /2 \sim \sqrt{1-t}$.
Therefore, estimate (\ref{est-Llam}) with $r=0$ is equivalent to
$$
\left| L_n^\lambda(t) \right| \le
     c_k \frac{n^{2 \lambda +1}}{(1+n \sqrt{1-t})^k}, \quad -1 \le t \le 1.
$$
Now, (\ref{est.main}) follows readily by Lemma~\ref{lem:loc1}.
\qed

\bigskip

The estimate of $|L_n^\mu(x,y)|$ from Theorem~\ref{thm:loc.ball}
enables us to control its $L^p$-norm.


\begin{prop} \label{prop:L-operator}
For $0< p \le \infty$, we have
\begin{equation} \label{est-norm-Ln}
\Big(\int_{B^d} |L_n^\mu(x,y)|^p W_\mu(y) dy\Big)^{1/p} \le c
\Big(\frac{n^{d}}{\W_\mu(n; x)}\Big)^{1-1/p},
\quad x\in B^d.
\end{equation}
\end{prop}

\begin{proof}
If $0 < p < \infty$ this proposition is an immediate consequence
of Theorem~\ref{thm:loc.ball} and Lemma~\ref{lem:est-norm} below,
taking into account that estimate (\ref{est.main}) holds for
an arbitrary~$k$.
In the case $p=\infty$ estimate (\ref{est-norm-Ln}) follows by
(\ref{est.main}) and (\ref{norm-dist2})
(see also estimate (\ref{Wx-Wy}) below).
\end{proof}


\begin{lem} \label{lem:est-norm}
If $\sigma > d/p + 2\mu|1/p-1/2|$, $\mu\ge 0$, $0 < p < \infty$, then
\begin{equation} \label{est-norm}
J_p:= \int_{B^d} \frac{W_\mu (y) dy}
   {\W_\mu(n;y)^{p/2} (1+n d(x,y))^{\sigma p}}
     \le c\,n^{-d}\W_\mu(n;x)^{1-p/2},
\end{equation}
where $c>0$ depends only on $p$, $\mu$, and $d$.
\end{lem}

\begin{proof} Let $\mu>0$ (the case $\mu=0$ is easier).
Three cases present themselves here.

{\em Case 1.} $p=2$.
Using spherical-polar coordinates and the fact that
$$
\int_{S^{d-1}}g(\langle x,y\rangle) d\omega(y) =
    \sigma_{d-2} \int_{-1}^1 g(|x| t) (1-t^2)^{(d-3)/2}dt,
$$
where $\sigma_{d-2}$ is the surface area of $S^{d-2}$, it follows that
\begin{align*}
J_2=
  c\int_0^1 \frac{r^{d-1} (1-r^2)^{\mu-1/2}}{(n^{-1}+\sqrt{1-r^2})^{2\mu}}
     \int_{-1}^1 \frac{(1-s^2)^{(d-3)/2} ds}{(1+n \arccos (r s|x| +
            \sqrt{1-|x|^2} \sqrt{1-r^2}\, ))^{2\sigma}} dr.
\end{align*}
Write briefly
$F(r,t) := 1/ [1+n \arccos (t |x| +
            \sqrt{1-|x|^2} \sqrt{1-r^2})]^{2\sigma}$.
Then
\begin{equation} \label{eq:Jmu}
  J_2 = c\int_0^1 \frac{r^{d-1} (1-r^2)^{\mu-1/2}}
  {(n^{-1} + \sqrt{1-r^2})^{2\mu}}
   \int_{-1}^1 F(r,rs)(1-s^2)^{(d-3)/2} ds dr.
\end{equation}
Next, we apply the substitution $u=rs$, then switch the order of integration,
and finally substitute $t=\sqrt{1-r^2}$.
This gives
\begin{align*}
J_2 &  = c\int_0^1 \frac{r (1-r^2)^{\mu-1/2}}{(n^{-1} + \sqrt{1-r^2})^{2\mu}}
   \int_{-r}^r F(r, u) (r^2-u^2)^{(d-3)/2} du dr \\
&  = c\int_{-1}^1 \int_{|u|}^1 F(r, u)
  \frac{r (1-r^2)^{\mu-1/2}}{(n^{-1} + \sqrt{1-r^2})^{2\mu}}
    (r^2-u^2)^{(d-3)/2} dr du \\
&  = c\int_{-1}^1 \int_{0}^{\sqrt{1-u^2}} F(\sqrt{1-t^2}, u)
  \frac{ t^{2\mu} (1-t^2-u^2)^{(d-3)/2}} {(n^{-1} + t)^{2\mu}}
      dt du.
\end{align*}
Using the trivial inequality $t /(t+n^{-1}) \le 1$
we conclude that
$$
J_2  \le c\int_{-1}^1 \int_{0}^{\sqrt{1-u^2}}
F(\sqrt{1-t^2}, u)(1-t^2-u^2)^{(d-3)/2} du dt.
$$
Since $\theta \sim \sin \theta/2 \sim \sqrt{1-\cos \theta}$ for
$0 \le \theta \le\pi$, we have
$$
F(\sqrt{1-t^2}, u) \sim
\left(1+n \sqrt{1- u |x|-t \sqrt{1-|x|^2}}\right)^{- 2\sigma},
\quad 0\le t\le \sqrt{1-u^2}.
$$
But
$1- u |x|-t \sqrt{1-|x|^2}\ge 0$ if $-\sqrt{1-u^2}\le t \le 0$.
Therefore, we can enlarge the domain of integration to obtain
$$
J_2  \le c\int_{B^2}
\frac{(1-t^2-u^2)^{(d-3)/2}du dt}
{\Big(1+n \sqrt{1- u |x|-t \sqrt{1-|x|^2}}\Big)^{2\sigma}}.
$$
Here $B^2$ is the unit disk in $\RR^2$.
We now change the variables $(u,t)\mapsto (a,b)$,
where
$$
  a = \sqrt{1-|x|^2}\,t  + |x| u, \qquad
  b = - |x| t + \sqrt{1-|x|^2}\,u.
$$
It is easy to see that this is an orthogonal transformation so that
$da\, db = du \,dt$.
Hence
\begin{align*}
J_2 &\le c\int_{B^2}\frac{(1-a^2 - b^2)^{(d-3)/2}} {(1+ n \sqrt{1-a})^{2\sigma}}
   da db\\
&= c\int_{-1}^1 \frac{1}{1+  n \sqrt{1-a})^{2\sigma}}
   \int_{-\sqrt{1-a^2}}^{\sqrt{1-a^2}} (1-a^2-b^2)^{(d-3)/2} db da \\
&\le c  \int_{-1}^1 \frac{(1-a^2)^{(d-2)/2}}
  {(1+n \sqrt{1-a})^{2\sigma}}   da\\
&\le \frac{c}{n^{2\sigma}} +
   c\int_{0}^1 \frac{(1-a)^{(d-2)/2}}{(1+n \sqrt{1-a})^{2\sigma}} da \\
&\le \frac{c}{n^{2\sigma}} +
    \frac{c}{n^d} \int_{0}^n \frac{s^{d-1}}{(1+s)^{2\sigma}} ds
   \le \frac{c}{n^d},
\end{align*}
since $2\sigma > d$. Thus (\ref{est-norm}) is established when $p=2$.

\smallskip
To prove (\ref{est-norm}) when $p \ne 2$ we will need the inequalities
\begin{align}\label{rel.norm-dist}
\frac{\sqrt{1-|x|^2}+n^{-1}}{\sqrt{2}(1+n\,d(x,y))}
&\le \sqrt{1-|y|^2}+n^{-1}\notag\\
&\le \sqrt{2}(\sqrt{1-|x|^2}+n^{-1})(1+n\,d(x,y)),
\quad x, y\in B^d,
\end{align}
which follow readily from (\ref{norm-dist2}).
From this and the definition of $\W_\mu(x;n)$ in (\ref{def.W}) we get
\begin{equation}\label{Wx-Wy}
c\W_\mu(n;x)(1+n\,d(x,y))^{-2\mu}
\le \W_\mu(n;y)
\le c\W_\mu(n;x)(1+n\,d(x,y))^{2\mu}.
\end{equation}


{\em Case 2.} $0<p<2$. Using (\ref{Wx-Wy}) we obtain
$$
\W_\mu(n;y)^{p/2} = \W_\mu(n;y)\W_\mu(n;y)^{p/2-1}
\ge \frac{c\W_\mu(n;y)}{\W_\mu(n;x)^{1-p/2}(1+nd(x,y))^{2\mu(1-p/2)}}
$$
and hence
$$
\int_{B^d} \frac{W_\mu(y)dy}{\W_\mu(n;y)^{p/2}(1+nd(x,y))^{\sigma p}}
\le c\W_\mu(n;x)^{1-p/2}
\int_{B^d} \frac{W_\mu(y)dy}{\W_\mu(n;y)(1+nd(x,y))^\tau},
$$
where $\tau :=(\sigma -2\mu(1/p-1/2))p$.
By the hypothesis of the lemma $\tau > d$. Then the above inequality and
(\ref{est-norm}) with $p=2$ imply (\ref{est-norm}) in this case.


\smallskip

{\em Case 3.} $2<p<\infty$.
Similarly as above by (\ref{Wx-Wy})
$$
\W_\mu(n;y)^{p/2} = \W_\mu(n;y)\W_\mu(n;y)^{p/2-1}
\ge \frac{c\W_\mu(n;y)\W_\mu(n;x)^{p/2-1}}{(1+nd(x,y))^{2\mu(p/2-1)}}.
$$
Consequently,
$$
\int_{B^d} \frac{W_\mu(y)dy}{\W_\mu(n;y)^{p/2}(1+nd(x,y))^{\sigma p}}
\le c\W_\mu(n;x)^{1-p/2}
\int_{B^d} \frac{W_\mu(y)dy}{\W_\mu(n;y)(1+nd(x,y))^\tau},
$$
where this time $\tau :=(\sigma -2\mu(1/2-1/p))p$.
Since $\tau > d$, the above inequality and
(\ref{est-norm}) with $p=2$ imply (\ref{est-norm}) in the case $2<p<\infty$.
\end{proof}

It will be vital for our further development that $L_n^\mu(x,y)$
is a $Lip \, 1$ function in $x$ (or $y$) with respect to the distance
$d(\cdot, \cdot)$. Throughout the rest of the paper,
we denote by $B_\xi (r)$ the closed ball centered at $\xi$ of radius $r >0$
with respect to the distance $d(\cdot,\cdot)$ on $B^d$, i.e.
$$
B_\xi(r):=\{x\in B^d: d(x, \xxi) \le r\}, \quad \xi \in B^d, \: r > 0.
$$


\begin{prop} \label{prop:Lip}
Let $\xxi, y \in B^d$. Then for all $x, z\in B_\xi(c^*n^{-1})$
$(c^*>0, n\ge 1)$ and an arbitrary $k$,
we have
\begin{equation} \label{Lip}
|L_n^\mu(x,y)-L_n^\mu(\xxi,y)|
\le c_k \frac{n^{d+1}d(x, \xxi)}
{\sqrt{\W_\mu(n; y)}\sqrt{\W_\mu(n; z)}(1+nd(y,z))^k},
\end{equation}
where $c_k$ depends only on $k$, $\mu$, $d$, $\ha$, and $c^*$.
\end{prop}

\begin{proof}Let $\mu>0$.
We will use the notation
$t(x,y;u):= \langle x, y\rangle + u\sqrt{1-|x|^2}\sqrt{1-|y|^2}$,
introduced in (\ref{def.t}).
By (\ref{rep.Ln}) it follows that
\begin{align}\label{difference}
&|L_n^\mu(x,y)-L_n^\mu(\xxi,y)|\notag\\
&\qquad\qquad\le c
\int_{-1}^1 \Big|L_n^\lambda(t(x,y;u))-L_n^\lambda(t(\xxi,y;u))\Big|
(1-u^2)^{\mu-1}du\\
&\qquad\qquad
\le c \int_{-1}^1
\Big\|
 \partial L_n^\lambda(\cdot)\Big\|_{L^\infty(I_u)}
|t(x,y;u)-t(\xxi,y;u)|(1-u^2)^{\mu-1}du,\notag
\end{align}
where $\partial f = f'$ and $I_u$ is the interval with end points
$t(x,y;u)$ and $t(\xxi,y;u)$.

As in the proof of Theorem~\ref{thm:loc.ball},
by estimate (\ref{est-Llam}) with $r=1$
it follows that
\begin{align}\label{est-deriv}
\Big\|\partial L_n^\lambda(\cdot)\Big\|_{L^\infty(I_u)}
&\le c_k n^{2\lambda+3} \max_{\tau\in I_u}\left (1+n\sqrt{1-\tau}\right)^{-k}
  \notag\\
&\le c_k n^{2\lambda+3}
\left( \left(1+n\sqrt{1-t(x,y;u)}\right)^{-k}+
  \left(1+n\sqrt{1-t(\xxi,y;u)}\right)^{-k}\right),
\end{align}
using the fact that $(1+n\sqrt{1-\tau})^{-k}$ is an increasing function
of $\tau$.

By the definition of $t(x,y;u)$ it follows that
(recall $x' := (x, \sqrt{1-|x|^2})$),
\begin{align*}
&|t(x,y;u)-t(\xxi,y;u)|\\
&\qquad\qquad\le |\langle x', y' \rangle - \langle \xxi', y' \rangle|
+ |1-u|\sqrt{1-|y|^2}\Big|\sqrt{1-|x|^2} - \sqrt{1-|\xxi|^2}\Big|\\
&\qquad\qquad\le |\cos d(x,y) - \cos d(\xxi,y)|
+ \sqrt{2}\,|1-u|\sqrt{1-|y|^2}d(x, \xxi),
\end{align*}
where we used inequality (\ref{norm-dist2}) from Lemma~\ref{lem:norm-dist}.
Denote briefly $\alpha:=d(x,y)$ and $\beta:=d(\xxi,y)$.
Then
\begin{align*}
|\cos d(x,y) - \cos d(\xxi,y)|
&=2\sin\frac{|\alpha-\beta|}{2}\sin\frac{\alpha+\beta}{2}
\le \frac{1}{2}|\alpha-\beta|(\alpha+\beta)\\
&\le \frac{1}{2}|d(x,y)-d(\xxi,y)|(d(x,y)+d(\xxi,y))\\
&\le \frac{1}{2}d(x,\xxi)(d(x,y)+d(\xxi,y))\\
&\le d(x,\xxi)(d(y,z)+c^*n^{-1})
\end{align*}
for $z\in B_\xi(c^*n^{-1})$. Hence
$$
|t(x,y;u)-t(\xxi,y;u)|\le d(x,\xxi)(d(y,z)+c^*n^{-1})+
\sqrt{2}\,|1-u|\sqrt{1-|y|^2}d(x, \xxi).
$$
We use this and (\ref{est-deriv}) in (\ref{difference}) to obtain
\begin{align*}
&|L_n^\mu(x,y)-L_n^\mu(\xxi,y)|
\le cn^{2\lambda+3}d(x,\xxi)(d(y,z)+c^*n^{-1})\\
& \quad\qquad \times
\Big(\int_{-1}^1 \frac{(1-u^2)^{\mu-1}du}{(1+n\sqrt{1-t(x,y;u})^k}
+\int_{-1}^1 \frac{(1-u^2)^{\mu-1}du}{(1+n\sqrt{1-t(\xxi,y;u})^k}
\Big)\\
& \qquad +cn^{2\lambda+3}\sqrt{1-|y|^2}d(x, \xxi)\\
& \quad\qquad \times\Big(
\int_{-1}^1 \frac{(1-u)(1-u^2)^{\mu-1}du}{(1+n\sqrt{1-t(x,y;u})^k}
+\int_{-1}^1 \frac{(1-u)(1-u^2)^{\mu-1}du}{(1+n\sqrt{1-t(\xxi,y;u})^k}
\Big)\\
& \qquad =:A_1+A_2+A_3+A_4.
\end{align*}
By Lemma~\ref{lem:loc1} with $\gamma=\mu-1$, we have
\begin{align*}
A_1\le cn^{2\lambda+3}d(x,\xxi)(d(y,z)+c^*n^{-1})
\frac{n^{-2\mu}}{\sqrt{\W_\mu(n;x)}\sqrt{\W_\mu(n;y)}(1+nd(x,y))^\sigma}
\end{align*}
with $\sigma:= k-3\mu-1$. Note that for $y\in B^d$ and all $z \in
B_\xxi(c^*n^{-1})$, we have
$1+nd(z,y) \sim 1+nd(\xxi, y)$ and
$\sqrt{1-|z|^2}+c^*n^{-1} \sim \sqrt{1-|\xxi|^2}+ n^{-1}$,
using (\ref{norm-dist2}).
Consequently,
\begin{align}\label{est.A1}
A_1\le \frac{cn^{d+1}d(x,\xxi)}
{\sqrt{\W_\mu(n;x)}\sqrt{\W_\mu(n;z)}(1+nd(y,z))^{\sigma-1}}.
\end{align}
We similarly obtain the same bound for $A_2$.

To estimate $A_3$ we employ Lemma~\ref{lem:loc1} with $\gamma=\mu$
and obtain
\begin{align}\label{est.A3}
&A_3\le cn^{2\lambda+3}\sqrt{1-|y|^2}d(x, \xxi)\notag\\
&\qquad\qquad\times \frac{n^{-2\mu-2}}
{(\sqrt{1-|x|^2}+n^{-1})^{\mu+1}(\sqrt{1-|y|^2}+n^{-1})^{\mu+1}
(1+nd(x,y))^{\sigma}}
\end{align}
with $\sigma:= k-3\mu-4$.
By cancelling appropriate terms we conclude that
(\ref{est.A1}) holds for $A_3$ as well. Exactly in the same way one can see
that $A_4$ also satisfies (\ref{est.A3}) and hence (\ref{est.A1}).
The proof of the proposition is complete.
\end{proof}

\noindent
{\bf Operators.}
We next use the localized polynomials $L_n^\mu$ as kernels
of linear operators defined by
\begin{equation} \label{def.Ln}
(\CL_n^\mu f)(x) := b_d^\mu \int_{B^d} f(y) L_n^\mu(x,y) W_\mu(y) dy,
\quad \mu \ge 0.
\end{equation}
Let $E_n^B(f)_{\mu,p}$ denote the best approximation
to $f \in L^p_\mu$, where $L^p_\mu:=L^p(B^d, W_\mu)$,
from the space $\Pi_n^d$ of all polynomials of degree
at most $n$, that is,
$$
E_n^B(f)_{\mu,p} := \inf_{q \in \Pi_n^d} \|f - q\|_{L^p_\mu}.
$$

\begin{prop} \label{thm:L-operator}
Let $\wh a$ be admissible of type $(a)$. Then the operator
$\CL_n^\mu$ satisfies the following properties:
\begin{enumerate}
\item[(i)] $\CL_n^\mu f$ is a polynomial of degree at most $2n$;
\item[(ii)] $\CL_n^\mu \,p = p$ for any polynomial $p$ of degree at most $n$;
\item[(iii)] for $f \in L^p_\mu$, $1 \le p \le \infty$,
\begin{equation}\label{est.operators}
\|\CL_n^\mu\|_{L^p_\mu\to L^p_\mu} \le c
\quad \hbox{and} \quad
\|\CL_n^\mu f - f \|_{L^p_\mu} \le c E_n^B(f)_{\mu,p}.
\end{equation}
\end{enumerate}
\end{prop}

\begin{proof}
The first two properties are obvious from the definition of $\CL_n^\mu$.
Since $\CL_n^\mu$ is an integral operator, the operator
norms $\|\CL_n^\mu\|_{L^1_\mu \to L^1_\mu}$ and
$\|\CL_n^\mu\|_{L^\infty \to L^\infty}$ are both bounded by
$$
\max_{x \in B^d} \int_{B^d}|L_n^\mu(x, y)|W_\mu(y)dy.
$$
Estimate (\ref{est-norm-Ln}) from Proposition \ref{prop:L-operator} with $p =1$
shows that this quantity is bounded by a constant
independent of $n$.
Then it follows by interpolation that
$\CL_n^\mu$ is a bounded operator from $L^p_\mu$ into $L^p_\mu$
for $1\le p \le \infty$, which yields (\ref{est.operators}).
\end{proof}

A result of the same nature holds true for more general weight functions of
the form $h_\kappa^2(x) (1-|x|^2)^{\mu-1/2}$, where $h_\kappa(x)$ is some
function invariant under a finite reflection group, see \cite{X04}.

\section{Cubature formula on $B^d$}\label{Cubature}
\setcounter{equation}{0}

Cubature formulae on $B^d$ with weights $W_\mu(x)$, $\mu \ge 0$, which are exact
for all polynomials of degree $n$ are valuable from many prospectives. Those
with positive coefficients are preferred for numerical computation and are called
positive cubature formulae. In the literature, only a handful of positive
cubature formulae are known.
For our purpose of constructing polynomial frames on $B^d$ (see \S\ref{Frame}))
we will need positive cubature whose knots are almost equally distributed
with respect to the distance $d(\cdot, \cdot)$ introduced in (\ref{def.dist}).
To the best of our knowledge there are no such
cubature formulae available up to now.
There is a close relation between cubature
formulae on the unit ball and those on the unit sphere $S^d$ \cite{X98}.
In the following we will follow the method used in \cite{MNW1} for constructing
cubature formulae on the unit sphere.

One of the difficulties in constructing cubature formulae on $B^d$ is the
lack of uniformly distributed points on $B^d$. We shall use as a substitute
sets of ``almost equally distributed points" with respect to the distance
$d(\cdot, \cdot)$ in $B^d$ which we describe in the following.



\begin{lem}\label{l:points}
For any $0<\eps \le \pi$ there exists
a partition $\cR_\eps$ of $B^d$
consisting of projections of spherical simplices
and a set $\cX_\eps \subset B^d$
{\rm(}consisting of their ``centers"{\rm)}
with the properties:
\begin{enumerate}
\item[(i)]
$B^d = \bigcup_{R \in \cR_\eps} R$ and the sets in $\cR_\eps$
do not overlap
$(R_1^\circ \cap R_2^\circ = \emptyset$ if $R_1 \ne R_2)$.
\item[(ii)]
For each $R \in \cR_\eps$ there is a unique $\xxi \in \cX_\eps$
such that
$
B_\xxi(c^*\eps) \subset R \subset B_\xxi(\eps).
$
\item[(iii)]
$\# \cX_\eps = \# \cR_\eps \le c^{**}\eps^{-d}$.
\end{enumerate}
Here $c^*$ and $c^{**}$ are constants depending only on $d$.
\end{lem}

\begin{proof}
As we already mentioned the distance
$d(x,y)$ ($x,y \in B^d$) is the geodesic distance
between $x', y'\in S^d_+$. So, we need to subdivide properly $S^d_+$.
We first divide $S^d_+$ into $2^{d}$ spherical simplices  analogous
to the intersections of $S^2$ with the octants in $\RR^3$.
Let $\cO_1$ be the spherical simplex on which all coordinates of
$\xi \in \cO_1$
are nonnegative and let
$$
\cT_1:=\Big\{\sum_{j=1}^{d+1}t_je_j: t_j \ge 0, \;
\sum_{j=1}^{d+1}t_j =1 \Big\},
$$
where $\{e_j\}$ are the standard unit vectors in $\RR^{d+1}$.
If $v:=(1, \dots, 1)$, then the map
$x(\xi):= \frac{\xi}{\langle \xi, v\rangle}$
gives an one-to-one correspondence between $\cO_1$ and $\cT_1$.
It~is readily seen that for $\xi, \zeta \in \cO_1$
\begin{equation}\label{relation-norms}
\frac{1}{2\sqrt{d}}d(\xi, \zeta) \le |x(\xi)-x(\zeta)|
\le 2\sqrt{d}\: d(\xi, \zeta).
\end{equation}
Here $|\cdot|$ denotes the Euclidean norm in $\RR^{d+1}$
and $d(\cdot,\cdot)$ is the geodesic distance on $S^d\subset \RR^{d+1}$.

We set $M:= \lceil2\sqrt{d}\eps^{-1}\rceil$ and divide the equilateral simplex
$\cT_1$ into $M^d$ equal equilateral subsimplices of side length
$L=\sqrt{2}/M$.
We denote by $\widetilde\cR_\eps^1$ the set of all spherical simplices obtained
by applying the inverse map $x^{-1}$ to the simplices defined above.
We similarly define the set $\widetilde\cX_\eps^1$ of the ``centers'' of all
spherical simplices
in $\widetilde\cR_\eps^1$ by applying the inverse map $x^{-1}$ to the midpoints
of the corresponding Euclidean simplices.
After these preparations,
we define $\cR_\eps^1$ as the set of projections onto $B^d$ of all spherical
simplices from $\widetilde\cR_\eps^1$ and we similarly define $\cX_\eps^1$.

It is straightforward to show that an equilateral Euclidean simplex
with each side of length $L$ contains the ball of radius
$L/\sqrt{2d(d+1)}$ centered at its midpoint and is contained in a ball
of radius $<L/\sqrt{2}$ with the same center.
Then (\ref{relation-norms}) yields that
the corresponding spherical simplex contains the spherical cap centered
at its center and of radius $L/(2d\sqrt{2(d+1)})$ and is contained in
a spherical cap with the same center and radius
$<\sqrt{2d}L \le 2\sqrt{d}/M\le \eps$.
This establishes property (ii) of Lemma~\ref{l:points} for the spherical
simplices in $\cR_\eps^1$.
Also, we have
$\# \cX_\eps^1 = \# \cR_\eps^1 = M^d \le (4\sqrt{d} \eps^{-1})^d$.

Repeating this procedure with all other initial simplices,
we establish the existence of the desired partition $\cR_\eps$.
\end{proof}


\begin{defn}\label{def.eq-dist-p}
A set $\cX_\eps \subset B^d$ which, along with an associated
partition $\cR_\eps$ of $B^d$, has the properties of the sets
$\cX_\eps$ and $\cR_\eps$ of Lemma~\ref{l:points} will be called a
set of {\em almost uniformly $\eps$-distributed points} on $B^d$.
\end{defn}

The following lemma contains additional information about the partition
$\cR_\varepsilon$.


\begin{lem}\label{lem:est.m-muR}
Let $\cR_\varepsilon$ be as in Lemma \ref{l:points}.
Then for any $\xi\in \cX_\eps$
\begin{equation}\label{est-mR}
|R_\xxi|:= \int_{R_\xxi} 1\,dx
\sim \eps^d\sqrt{1-|\xxi|^2}
\end{equation}
and
\begin{equation}\label{est-m-muR}
m_\mu(R_\xxi):= \int_{R_\xxi} W_\mu(x)\,dx
\sim \eps^d(1-|\xxi|^2)^\mu
= \eps^d\frac{W_\mu(\xxi)}{W_0(\xxi)}
\sim \eps^d(\sqrt{1-|\xxi|^2}+\eps)^{2\mu}.
\end{equation}
Here the constants of equivalence depend only on $d$ and $\mu$.
\end{lem}

\begin{proof}
To prove (\ref{est-mR}) we use property (ii) in Lemma \ref{l:points}
which yields
\begin{equation}\label{Rxi-Bxi}
|R_\xxi| \sim |B_\xxi(\eps)|
\quad \mbox{and} \quad
d(\xxi, \partial B^d)\ge c^*\eps.
\end{equation}
We can assume without
loss of generality that $\xxi$ lies on the positive $x_1$-axis, i.e.
$\xxi =(\xxi_1, 0, \dots, 0)$ and $0 < \xxi_1 < 1$. The boundary
$\partial B_\xxi(\eps)$ of $B_\xxi(\eps)$ is given by the equation
$x_1 \xi_1 + \sqrt{1-|x|^2} \sqrt{1-\xi_1^2} = \cos \varepsilon$.
A simple manipulation of this identity shows that
$\partial B_\xxi(\eps)$ is the ellipsoid
$$
\frac{(x_1-\xxi_1 \cos \varepsilon)^2}{1-|\xxi|^2}+
x_2^2+\dots +x_d^2=\sin^2 \varepsilon.
$$
From this it follows that $|B_\xxi(\eps)|\sim
\eps^d\sqrt{1-|\xxi|^2}$
(using that $\sin\eps \sim\eps$) and then (\ref{est-mR}) follows.


We now turn to the proof of (\ref{est-m-muR}).
There are two cases to be considered.

{\em Case 1.} $\mu \ge 1/2$.
Denote $R_\xxi^-:= R_\xxi\cap \{x\in B^d: |x|\le|\xxi|\}$.
It is easily seen that $|R_\xxi^-|\sim  |R_\xxi| \sim
\eps^d\sqrt{1-|\xxi|^2}$. Then
$$
 \int_{R_\xi} W_\mu(x)\,dx \ge  \int_{R_\xi^-} W_\mu(x)\,dx \ge
   W_\mu(\xxi)|R_\xxi^-| \sim W_\mu(\xxi)\eps^d\sqrt{1-|\xxi|^2}
= \eps^d (1-|\xxi|^2)^\mu.
$$
Since $\xi$ is in the center of $\cR_\xi$ by construction, we have
$\sqrt{1-|\xi|^2} \ge c \varepsilon$. Hence, for $x \in R_\xi
\subset B_\xi(\eps)$, inequality \eqref{norm-dist2} shows that
$$
W_\mu(x) \le (\sqrt{1-|\xi|^2}+ \varepsilon)^{2\mu-1} \le c W_\mu(\xi),
$$
which yields
$$
 \int_{R_\xi} W_\mu(x)\,dx \le cW_\mu(\xxi)|R_\xi|
   \sim W_\mu(\xxi)\eps^d\sqrt{1-|\xxi|^2}
    = \eps^d (1-|\xxi|^2)^\mu.
$$

{\em Case 2.} $0\le \mu < 1/2$.
Denote
$R_\xxi^+:= R_\xxi\cap \{x\in B^d: |x|\ge|\xxi|\}$.
Proceeding as above we again get (\ref{est-m-muR}).

Finally, using (\ref{Rxi-Bxi}) we obtain
$\sqrt{1-|\xxi|^2}\ge \sin c^*\eps \ge c\eps$
which implies the last equivalence in (\ref{est-m-muR}).
The proof of the lemma is complete.
\end{proof}


\begin{thm}\label{thm:cubature}
There exists a constant $c^\diamond >0$ {\rm(}depending only on $d${\rm)}
such that for any $n\ge 1$ and
a set $\cX_\eps$ of almost uniformly $\eps$-distributed points on $B^d$
with $\eps:=c^\diamond/n$,
there exist positive coefficients $\{\lambda_\xxi\}_{\xxi \in \cX_\eps}$
such that the cubature formula
$$
\int_{B^d} f(x)\,dx \sim \sum_{\xxi \in \cX_\eps} \lambda_\xxi f(\xxi)
$$
is exact for all polynomials of degree $\le n$. In addition,
$$
\lambda_\xxi \sim n^{-d}\W_\mu(n;\xi)\sim \eps^d (1-|\xxi|^2)^\mu
\sim m_\mu(B_\xxi(\eps))
$$
with constants of equivalence depending
only on $\mu$ and $d$. Here $m_\mu(E):=\int_E W_\mu(x)dx$.
\end{thm}

Note that when $\mu=0$ the cubature formula of Theorem~\ref{thm:cubature}
can be derived from the cubature formula on $S^{d+1}$ from \cite{MNW1, MNW2, NPW}
by applying \cite[Theorem 4.2]{X98}.

For the proof Theorem~\ref{thm:cubature} we will utilize the idea used
in \cite{MNW1, MNW2} (see also \cite{NPW})
for the construction of a cubature formula on $S^d$.

Assume that $\cX_\eps$ (with associated partition $\cR_\eps$)
is a set of almost uniformly $\eps$-distributed points on $B^d$
(see Definition~\ref{def.eq-dist-p}),
where $\eps = \delta/n$ with $n\ge 1$ and $\delta$
will be selected later on.
We introduce the following weighted
$\ell_1$-norm for functions defined on $B^d$:
\begin{equation}\label{weithed-norm}
\|f\|_{\ell^1_\mu(\cX_\eps)}:= \sum_{\xxi \in \cX_\eps}
|f(\xxi)|m_\mu(R_\xxi).
\end{equation}
Also, recall the notation
$
\|f\|_{L^1_\mu}=
\|f\|_{L^1(W_\mu, B^d)}:= \int_{B^d} |f(x)|W_\mu(x)dx.
$

We need a couple of additional results.


\begin{lem}\label{l:ell-norm}
If $g \in \Pi_n^d$, then
\begin{equation}\label{ell-norm1}
\left|\|g\|_{L^1_\mu}-\|g\|_{\ell^1_\mu(\cX_\eps)}\right|
\le \sum_{\xxi\in\cX_\eps}\int_{\R_\xxi}|g(x)-g(\xxi)|W_\mu(x)dx
\le c^\star\delta\|g\|_{L^1_\mu}
\end{equation}
and hence
\begin{equation}\label{ell-norm2}
(1-c^\star\delta)\|g\|_{L^1_\mu} \le \|g\|_{\ell^1_\mu(\cX_\eps)}
\le (1+c^\star\delta)\|g\|_{L^1_\mu},
\end{equation}
where $c^\star$ depends only on $d$ and $\mu$.
\end{lem}

\begin{proof}
Let $\CL_n^\mu$ be the operator defined in (\ref{def.Ln}).
By Proposition \ref{prop:L-operator} we have $g= \CL_n^\mu g$.
Using this and the fact that $\cR_\eps$ is a partition of $B^d$
(see Lemma~\ref{l:points}), we obtain
\begin{align*}
& \left|\|g\|_{L^1_\mu}-\|g\|_{\ell^1_\mu(\cX_\eps)}\right|
\le \sum_{\xxi\in\cX_\eps}\int_{\R_\xxi}|g(x)-g(\xxi)|W_\mu(x)dx\\
&\qquad\qquad
\le \sum_{\xxi\in\cX_\eps}\int_{\R_\xxi}\int_{B^d}
|L_n^\mu(x,y) - L_n^\mu(\xxi, y)||g(y)|W_\mu(y)dyW_\mu(x)dx\\
&\qquad\qquad
\le \|g\|_{L^1_\mu}\sup_{y\in B^d}
\sum_{\xxi\in\cX_\eps}\int_{\R_\xxi}
|L_n^\mu(x, y) - L_n^\mu(\xxi, y)|W_\mu(x)dx.
\end{align*}
By Proposition~\ref{prop:Lip} with $z=x$, it follows that
\begin{align*}
&\int_{\R_\xxi}|L_n^\mu(x, y) - L_n^\mu(\xxi, y)|W_\mu(x)dx \\
& \qquad\qquad
 \le \int_{\R_\xxi}
\frac{c_kn^{d+1}d(x, \xxi)W_\mu(x)dx}
{\sqrt{\W_\mu(n;x)}\sqrt{\W_\mu(n;y)}(1+nd(y,x))^k}.
\end{align*}
Choosing $k$ sufficiently large ($k>d+\mu$ will do)
we apply Lemma~\ref{lem:est-norm} with $p=1$ and use that
$d(x, \xxi) \le \delta n^{-1}$ for $x\in \R_\xxi$
to obtain
\begin{align*}
&\sup_{y\in B^d}
\sum_{\xxi\in\cX_\eps}\int_{\R_\xxi}
|L_n^\mu(x, y) - L_n^\mu(\xxi, y)|W_\mu(x)dx\\
&\qquad\qquad
\le c\delta n^d\int_{B^d} \frac{W_\mu(x)dx}
{\sqrt{\W_\mu(n;x)}\sqrt{\W_\mu(n;y)}(1+nd(y,x))^k}
\le c\delta.
\end{align*}
The lemma follows.
\end{proof}

\noindent
{\bf The Farkas Lemma.}
A variant of the well known in Optimization Farkas lemma
will play an important role in the proof of Theorem~\ref{thm:cubature}.

\begin{prop}\label{prop:Farkas}
Let $V$ be a finite dimensional real vector space and denote
by $V^*$ its dual.
Let
$u_1, u_2, \dots, u_n \in V^*$
and suppose
$u\in V^*$ has the property that
$u(x)\ge 0$ for all $x\in V$ such that
$u_j(x)\ge 0$ for $j=1, 2, \dots, n$.
Then there exist $a_j\ge 0$, $j=1, 2, \dots, n$, such that
\begin{equation}\label{Farkas}
u=\sum_{j=1}^n a_j u_j.
\end{equation}
\end{prop}
For the proof of this proposition, see e.g. \cite{Bor}.


\medskip

\noindent
{\it Proof of Theorem \ref{thm:cubature}.}
First, we choose $\delta~:=~\frac{1}{3c^\star}$,
where $c^\star$ is the constant from
Lemma~\ref{l:ell-norm}.
In applying Proposition~\ref{prop:Farkas}, we take
$V:= \Pi_n^d$ and
$\{u_j\}$ to be the set of all point evaluation functionals
$\{\delta_\xxi\}_{\xxi \in \cX_\eps}$.

Let the linear functionals $u$ and $u_\gamma$ be defined by
$$
u(g):= \int_{B^d} g(x)W_\mu(x)\,dx
\quad \mbox{and} \quad
u_\gamma(g):= u(g) - \gamma\sum_{\xxi \in \cX_\eps} g(\xxi)m_\mu(R_\xxi).
$$
Since $c^\star\delta=1/3$,
the left-hand-side estimate in (\ref{ell-norm2}) yields
\begin{equation}\label{L1<ell}
\|g\|_{L^1_\mu} \le (3/2)\|g\|_{\ell^1_\mu(\cX_\eps)},
\quad g \in \Pi_n^d.
\end{equation}
Suppose $g \in \Pi_n^d$ and $g(\xxi) \ge 0$ for all $\xxi \in \cX_\eps$.
Then using (\ref{ell-norm1}) with $c^*\delta=1/3$ and (\ref{L1<ell}),
we obtain
\begin{eqnarray*}
\left|u(g)-\|g\|_{\ell^1_\mu(\cX_\eps)}\right|
= \sum_{\xxi\in\cX_\eps}\int_{\R_\xxi}|g(x)-g(\xxi)|W_\mu(x)\,dx
\le c^\star\delta \|g\|_{L^1_\mu}
\le (1/2)\|g\|_{\ell^1_\mu(\cX_\eps)}
\end{eqnarray*}
and hence
$
u(g) \ge (1/2)\|g\|_{\ell^1_\mu(\cX_\eps)}.
$
Choose $\gamma:= 1/3$. Then since
$g(\xxi) \ge 0$, $\xxi \in \cX_\eps$,
we obtain
$$
u_\gamma(g) = u(g)-(1/3)\|g\|_{\ell^1_\mu(\cX_\eps)}
\ge (1/6)\|g\|_{\ell^1_\mu(\cX_\eps)} \ge 0.
$$
Applying Proposition \ref{prop:Farkas} to $u_\gamma$,
there exist numbers $a_\xxi \ge 0$, $\xxi \in \cX_\eps$, such that
$$
u_\gamma(g) = \sum_{\xxi \in \cX_\eps} a_\xxi g(\xxi),
\quad g \in \Pi_n^d,
$$
and hence
$$ 
u(g) = \sum_{\xxi \in \cX_\eps} (a_\xxi + (1/3)m_\mu(R_\xxi)) g(\xxi)
=: \sum_{\xxi \in \cX_\eps} \lambda_\xxi g(\xxi), \quad
g \in \Pi_n^d.
$$ 
Therefore, the linear functional
$\sum_{\xxi \in \cX_\eps} \lambda_\xxi g(\xxi)$
provides a cubature formula exact for all polynomials of degree $n$.

Clearly, $\lambda_\xxi \ge m_\mu(R_\xxi)/3$ and the estimate
$\lambda_\xxi\le cm_\mu(R_\xxi)$ follows from
Lemma~\ref{lem:est.m-muR} and Proposition \ref{prop:coeff} below.
\qed

\medskip


The last ingredient in bounding $\lambda_\xxi$ from above is the following
general result that is of independent interest.


\begin{prop}\label{prop:coeff}
If a positive cubature formula
\begin{equation}\label{positive-cubature}
\int_{B^d} f(x) W_\mu(x) dx
\sim \sum_{\xi \in \cX_\varepsilon} \lambda_\xi f(\xi),
\quad \lambda_\xi >0, ~~~|\xi| < 1,
\end{equation}
is exact for all polynomials of degree $\le n$, then
\begin{equation}\label{est-lam-xi}
\lambda_\xi \le  cn^{-d}\W_\mu(n;\xi)
=cn^{-d}(\sqrt{1-|\xi|^2} + n^{-1})^{2\mu},
\quad \xi\in \cX_\eps,
\end{equation}
where $c>0$ depends only on $\mu$ and $d$.
\end{prop}

\begin{proof}
Recall the kernel $K_m(W_\mu;x,y)$ defined in \eqref{def.Kn}. Evidently
$K_m(W_\mu;\xi,\xi) > 0$ and
$$
  \int_{B^d} \left[K_m(W_\mu;x,y)\right]^2 W_\mu(y)dy = K_m(W_\mu;x,x).
$$
Let $m = \lfloor n/2 \rfloor$. Then it follows that
$$
\lambda_\xi  \le \sum_{\eta \in \cX_\varepsilon} \lambda_\eta
   \left [\frac{K_m(W_\mu;\xi,\eta)}{K_m(W_\mu;\xi,\xi)}\right]^2
  = \int_{B^d} \left[\frac{K_m(W_\mu;\xi,x)}{K_m(W_\mu;\xi,\xi)}\right]^2
   W_\mu(x) dx
  = \frac{1}{K_m(W_\mu; \xi,\xi)}.
$$
Hence, the stated result is a consequence of an upper bound
for $[K_m(W_\mu;x,x)]^{-1}$, to be established in
Proposition \ref{prop:Christoffel} below.
\end{proof}

In order to establish the needed upper bound for $[K_n(W_\mu;x,x)]^{-1}$
we now construct a family of well localized polynomials.


\begin{lem}\label{lem:local.polyn}
For any $k, m\ge1$ and $\xi\in B^d$ there exist a polynomial
$P_\xi\in \Pi^d_{2km}$ and a constant $c^*>0$ depending only on $k$ and $d$
such that $P_\xi(\xi)=1$ and for $0\le\gamma\le k$, $x\in B^d$,
\begin{equation}\label{est.Pxi}
0\le P_\xi(x)\le \frac{c^*}{(1+md(\xi, x))^{2k}}
\le \frac{c(\sqrt{1-|\xi|^2} + m^{-1})^\gamma}
            {(\sqrt{1-|x|^2} + m^{-1})^\gamma(1+md(\xi, x))^k}.
\end{equation}
\end{lem}

\begin{proof}
Let $q(\theta):= \Big(\frac{\sin (m\theta/2)}{m\sin (\theta/2)}\Big)^{2k}$.
Evidently, $q$ is a trigonometric polynomial of degree less than $km$,
$q(0)=1$, and
\begin{equation}\label{ets.q}
0\le q(\theta)\le \frac{c}{(1+m|\theta|)^{2k}}, \quad |\theta|\le \pi.
\end{equation}
For $0\le \alpha\le \pi$, we define the algebraic polynomial
$Q_\alpha(t)$ by
$$
Q_\alpha(\cos\theta):=\frac{q(\theta-\alpha)+q(\theta+\alpha)}{1+q(2\alpha)}.
$$
It is readily seen that
$\deg Q_\alpha < km$, $Q_\alpha(\cos \alpha)=1$, and
\begin{equation}\label{ets.Qalpha}
0\le Q_\alpha(\cos\theta)
\le \frac{c}{(1+m|\theta-\alpha|)^{2k}}, \quad 0\le \theta \le \pi.
\end{equation}
Also, $Q_{\pi/2}$ is even and
\begin{equation}\label{ets.Qpi}
0\le Q_{\pi/2}(t)
\le \frac{c}{(1+m|\arccos t-\pi/2|)^{2k}}
\le \frac{c}{(1+m|t|)^{2k}}, \quad |t|\le 1.
\end{equation}
Without loss of generality we may assume that
$\xi=(\xi_1, 0, \dots, 0)$ with $0 < \xi_1 < 1$.
We choose $\alpha \in (0, \pi/2)$ so that $\xi_1=\cos\alpha$.
Then (\ref{ets.Qalpha}) gives
\begin{equation}\label{ets.Qalpha1}
0\le Q_\alpha(t)
\le \frac{c}{(1+m\,d_1(\xi_1, t))^{2k}}, \quad |t|\le 1,
\end{equation}
where
$d_1(\xi_1, t):=\arccos\,(\xi_1t+\sqrt{1-\xi_1^2}\sqrt{1-t^2})$
is the univariate version of the distance $d(\cdot, \cdot)$
(see (\ref{def.dist})).
We define
$$P_\xi(x):=Q_\alpha(x_1)Q_{\pi/2}\Big(\sqrt{x_2+\cdots+x_d^2}\Big).$$
Clearly, $P_\xi \in \Pi^d_{2km}$, $P_\xi(\xi)=1$,
and by (\ref{ets.Qpi})-(\ref{ets.Qalpha1})
\begin{equation}\label{est.Pxi1}
0\le P_\xi(x) \le \frac{c}{[(1+m|x_*|)(1+md_1(\xi_1,x_1))]^{2k}},
\quad x\in B^d,
\end{equation}
where $x_*:=(x_2, \dots, x_d)$ and
$|x_*|:=(x_2^2+\cdots+x_d^2)^{1/2}$.

It remains to show that $P_\xi$ obeys (\ref{est.Pxi}).
To this end we fist show that
\begin{equation}\label{d<norm+d1}
d(\xi, x)\le 2(|x_*|+d_1(\xi_1, x_1)).
\end{equation}
Denote briefly $x_\diamond:=(x_1, 0, \dots, 0)$.
We have
$$
d(\xi, x)\le d(\xi, x_\diamond)+ d(x_\diamond, x)
= d_1(\xi_1, x_1)+d(x_\diamond, x).
$$
Our next step is to prove the inequality
\begin{equation}\label{d<x*}
d(x_\diamond, x)\le 2|x_*|.
\end{equation}
Evidently,
$$
d(x_\diamond, x)=\arccos(\langle x_\diamond', x'\rangle)
= \arccos \Big(x_1^2+\sqrt{1-x_1^2}\sqrt{1-x_1^2-x_2^2-\cdots-x_d^2}\Big).
$$
One easily verifies the inequality
$\arccos t \le 2\sqrt{1-t}$, $0\le t \le 1$,
and hence (\ref{d<x*}) will be established if we show that
$$
\Big(1-x_1^2-\sqrt{1-x_1^2}\sqrt{1-x_1^2-|x_*|^2}\Big)^{1/2}
\le |x_*|.
$$
Denote briefly $a:=\sqrt{1-x_1^2}$ and $b:=|x_*|$.
Then the above inequality is equivalent to
$a^2-a\sqrt{a^2-b^2} \le b^2$ or $a\sqrt{a^2-b^2}-(a^2-b^2)\ge 0$.
But the latter inequality is apparently valid since
$$
a\sqrt{a^2-b^2}-(a^2-b^2)=\frac{b^2\sqrt{a^2-b^2}}{a+\sqrt{a^2-b^2}}\ge 0.
$$
Thus (\ref{d<x*}) is established and hence (\ref{d<norm+d1}) holds.
Combining (\ref{est.Pxi1}) with (\ref{d<norm+d1}) gives
\begin{equation}\label{est.Pxi2}
0\le P_\xi(x) \le \frac{c}{(1+m\,d(\xi,x))^{2k}},
\quad x\in B^d,
\end{equation}
which is the first estimate of $P_\xi(x)$ in (\ref{est.Pxi}).

To prove the second estimate in (\ref{est.Pxi}) we need the estimate
\begin{equation}\label{est.dL}
\frac{1}{1+m\,d(\xi,x)}
\le c\frac{\sqrt{1-|\xi|^2}+m^{-1}}{\sqrt{1-|x|^2}+m^{-1}},
\quad x\in B^d,
\end{equation}
which apparently follows by inequality (\ref{norm-dist2}) in
Lemma~\ref{lem:norm-dist} (see also (\ref{rel.norm-dist})).

Finally, applying (\ref{est.dL}) in (\ref{est.Pxi2}) we get the second
estimate in (\ref{est.Pxi}),
which completes the proof.
\end{proof}

The function $\Lambda_n(x):= [K_n(W_\mu;x,x)]^{-1}$ is the so called
Christoffel function, which has the following characteristic
property \cite[p. 109]{DX}:
\begin{equation}\label{Christoffel}
\Lambda_n(x) = \min_{P(x) = 1, P\in \Pi_n^d} \int_{B^d}[P(y)]^2 W_\mu(y) dy,
\quad x \in B^d.
\end{equation}
The localized polynomials in Lemma~\ref{lem:local.polyn} give an upper bound
for the Christoffel function, used in the proof of Proposition
\ref{prop:coeff}.


\begin{prop}\label{prop:Christoffel}
For any $\mu \ge 0$ and $d>1$ there exists a constant $c>0$
such that
\begin{equation}\label{Christoffel1}
\Lambda_n(x)\le cn^{-d}\W_\mu(n;x), \quad x\in B^d, \: n\ge 1.
\end{equation}
\end{prop}

\begin{proof}
Write $k:=[\max\{d/2, \mu\}]+1$ and let $n\ge 4k$
(the case $1 \le n< 4k$ is trivial).
Set $m:=[n/2k]$.
By Lemma~\ref{lem:local.polyn} there exists a polynomial
$P_x(y) \in \Pi^d_n$ such that $P_x(x)=1$ and (\ref{est.Pxi}) holds
with $\gamma=\mu$ and $\xxi, x$ replaced by $x, y$.
Then by (\ref{Christoffel}), (\ref{est.Pxi}),
and Lemma~\ref{lem:est-norm} with $p=2$,
we infer
\begin{align*}
\Lambda_n(x)&\le \int_{B^d}[P_x(y)]^2W_\mu(y)\,dy
\le c\int_{B^d}\frac{\W_\mu(m;x)W_\mu(y)dy}{\W_\mu(m;y)(1+m\,d(x, y))^{2k}}\\
&\le cm^{-d}\W_\mu(m;x)\le cn^{-d}\W_\mu(n;x).
\end{align*}
\end{proof}

For the construction of our frames,
we will need the following result which is
an immediate consequence of Lemma~\ref{l:points} and
Theorem~\ref{thm:cubature}.


\begin{cor}\label{c:cubature}
There exists a sequence $\{\cX_j\}_{j=0}^\infty$
of sets of almost uniformly $\eps_j$-distributed points on $B^d$
$(\cX_j:= \cX_{\eps_j})$ with $\eps_j:= c^\diamond 2^{-j-2}$
and there exist positive coefficients
$\{\lambda_\xxi\}_{\xxi \in \cX_j}$
such that the cubature
\begin{equation}\label{cubature-j}
\int_{B^d} f(x)W_\mu(x)\,dx \sim \sum_{\xxi \in \cX_j} \lambda_\xxi f(\xxi)
\end{equation}
is exact for all polynomials of degree $\le 2^{j+2}$.
Moreover, $\lambda_\xxi \sim m_\mu(B_\xxi(2^{-j}))$ and
$\# \cX_j \sim 2^{jd}$
with constants of equivalence depending only on $d$ and $\mu$.
\end{cor}

\section{Tight polynomial frames (needlets) in $L^2(B^d, W_\mu)$}
\label{Frame}
\setcounter{equation}{0}

We will apply the general scheme, described in \S\ref{Principles},
for construction of polynomial frames in $L_\mu^2:=L^2(B^d, W_\mu)$.
To this end we will utilize the
localized polynomials from Theorem~\ref{thm:loc.ball}
and the cubature formula from Theorem~\ref{thm:cubature}
(see Corollary~\ref{c:cubature}).

Let $\wh a$ satisfy the conditions:
\begin{equation}\label{a1}
\wh a \in C^\infty[0, \infty), \quad \wh a \ge 0,
 \quad \supp \wh a \subset [1/2, 2],
\end{equation}
\begin{equation}\label{a2}
\wh a(t)>c>0, \qquad \mbox{if $t \in [3/5, 5/3]$},
\end{equation}
\begin{equation}\label{a3}
\wh a^2(t) + \wh a^2(2t) =1,
\qquad \mbox{if $t \in [1/2, 1]$.}
\end{equation}
For the construction of such functions, see e.g. \cite{NPW}.

We introduce a sequence of polynomial ``kernels'':
$L_0(x,y):=1$ and (see \S\ref{Loc.pol.ball})
$$ 
L_j(x,y):= \sum_{\nu=0}^\infty
\wh a \Big(\frac{\nu}{2^{j-1}}\Big)P_\nu(W_\mu;x,y),
\quad  j=1, 2, \dots.
$$ 
We now define the needlets (frame elements) by
$$ 
\psi_\xi(x):= \sqrt{\lambda_\xi}\cdot L_j(x, \xi)
\quad \mbox{for}\quad
\xi\in \cX_j, ~ j=0, 1, \dots,
$$ 
where $\cX_j$ is the set of the knots and the $\lambda_\xi$'s are
the coefficients
of the cubature formula \eqref{cubature-j} from Corollary~\ref{c:cubature}.
We write
$\cX := \cup_{j = 0}^\infty \cX_j$ (see \S\ref{Principles})
and define the needlet system $\Psi$ by
$$
\Psi:=\{\psi_\xi\}_{\xi\in\cX}.
$$

Denoting
$$
(L_j*f)(x):= \int_{B^d}L_j(x,y)f(y)W_\mu(y)dy,
$$
we get as in (\ref{L.reprG}) the semi-discrete needlet decomposition
of $L^2_\mu$:
$$ 
f=\sum_{j=0}^\infty L_j*L_j*f
\quad\mbox{for $f\in L^2_\mu$.}
$$ 

The next theorem shows that the needlet system
$\Psi$ is a tight frame in $L^2_\mu$.


\begin{thm}\label{t:frame}
If $f\in L^2_\mu$, then
$$ 
f =\sum_{j=0}^\infty\sum_{\xi\in\cX_j}
  \langle f, \psi_{\xi} \rangle \psi_{\xi}
= \sum_{\xi\in \cX} \langle f, \psi_\xi\rangle \psi_{\xi}
\quad\mbox{in $L^2_\mu$}
$$ 
and
$$ 
\|f\|_{L^2_\mu}
=\Big(\sum_{\xi\in \cX} |\langle f, \psi_\xi\rangle|^2\Big)^{1/2}.
$$ 
\end{thm}
This theorem follows at once from Proposition~\ref{prop:frameG}
and Theorem~\ref{thm:cubature} (see Corollary~\ref{c:cubature}).

We finally show that each needlet $\psi_\xi$
has faster than any polynomial rate of decay
away from its center (pole) $\xi$ with respect to the distance
$d(\cdot, \cdot)$ on $B^d$, which prompted us to coin their name.
This property of the needlets is critical for using them
for decomposition of spaces other than $L_\mu^2$.


\begin{thm}\label{p:localization}
For any $k>0$ there exists a constant $c_k>0$ depending only on
$k$, $\mu$, $d$, and $\wh a$ such that for $\xi\in \cX_j$, $j=0, 1, \dots$,
\begin{equation} \label{est.needl}
|\psi_\xi(x)| \le c_k \frac{2^{jd/2}}{\sqrt{\W_\mu(2^j;x)}
(1 + 2^jd(x,\xi))^k}, \quad x\in B^d.
\end{equation}
\end{thm}

\begin{proof}
Estimates (\ref{est.needl}) follows readily from (\ref{est.main})
(see Theorem~\ref{thm:loc.ball}), taking into account that
$\lambda_\xi \le c2^{-jd}\W_\mu(2^j;\xi)$
for $\xi \in \cX_j$.
\end{proof}


\begin{re}
Estimate $(\ref{est.needl})$
and Lemma~\ref{lem:est-norm} with $p = 2$
yield
$\|\psi_{\xi}\|_{L^2_\mu}\le c$
for all $\xi\in \cX$,
which shows that estimate $(\ref{est.needl})$ is sharp (in a sense).
\end{re}

\section{Appendix}\label{appendix}
\setcounter{equation}{0}


\noindent
{\em Proof of Theorem~\ref{thm:BD} $($see \cite{BD}$)$.}
We will need the following well known estimate for Jacobi polynomials
\cite[(7.32.6), p. 167]{Sz}: For $\a, \b >-1/2$ and $n\ge 1$,
\begin{equation}\label{est.Pn}
|P_n^{(\a,\b)}(\cos \theta)| \le
c(\a,\b)\left\{\aligned
&\min\{n^\a, n^{-1/2}\theta^{-\a-1/2}\} &\mbox{if $0\le \theta\le \pi/2$,}\\
&\min\{n^\b, n^{-1/2}(\pi-\theta)^{-\b-1/2}\} &\mbox{if $\pi/2\le \theta\le \pi$.}
\endaligned
\right.
\end{equation}

We first prove (\ref{est.Ln}) for $r=0$.
We may assume that $n >2k$ and $k\ge 2$, since (\ref{est.Ln}) is trivial
when $n \le 2k$ ($c_k$ may depend on $k$).

{\em Case 1.} $0 \le \theta \le \pi/n$.
Since $\Gamma(j + a) /\Gamma(j+1) \sim j^{a-1}$ 
and $\a\ge\b>-1/2$, it follows by (\ref{est.Pn}) that
\begin{align*}
|L_n^{\a,\b}(\cos\theta)|
     \le c \sum_{j=0}^{2n} j^{2 \a +1} \le c n^{2 \a+2},
\end{align*}
which yields (\ref{est.Ln}) for $0\le \theta\le \pi/n$.

{\em Case 2.} $\pi/n \le \theta \le \pi$.
A key role here will be played by the identity \cite[(4.5.3), p.71]{Sz}:
\begin{align}\label{key-ident}
&\sum_{\nu=0}^n \frac{(2\nu+\a+k+\b+1)\Gamma(\nu+\a+k+\b+1)}{\Gamma(\nu+\b+1)}
P_\nu^{(\a+k,\b)}(x)\\
&\qquad\qquad= \frac{\Gamma(n+\a+k+1+\b+1)}{\Gamma(n+\b+1)}
P_n^{(\a+k+1,\b)}(x).\notag
\end{align}
Applying summation by parts to the sum in (\ref{eq:L-n})
(using (\ref{key-ident}) with $k=0$), we get
\begin{equation}\label{Ln-parts}
L_n^{\a,\b}(x) =  c^\diamond
\sum_{j=0}^\infty
\Big[\wh a\Big(\frac{j}{n}\Big)-\wh a\Big(\frac{j+1}{n}\Big)\Big]
\frac{\Gamma(j+\a+1+\b+1)}{\Gamma(j+\b+1)}P_j^{(\a+1,\b)}(x).
\end{equation}
We now define the sequence of functions $(A_k(t))_{k=0}^\infty$ by
$$
A_0(t):= (2t+\a+\b+1)\wh a\Big(\frac{t}{n}\Big)
$$
and inductively
\begin{equation}\label{Ak+1}
A_{k+1}(t):= \frac{A_k(t)}{2t+\a+k+\b+1}-\frac{A_k(t+1)}{2t+\a+k+\b+3},
\quad k\ge 0.
\end{equation}
It is readily seen that
\begin{equation}\label{A1}
A_1(t):= \wh a\Big(\frac{t}{n}\Big)-\wh a\Big(\frac{t+1}{n}\Big)
\end{equation}
and hence
$\supp A_k \subset [n-k, 2n] \subset [n/2, 2n]$, $1\le k \le n/2$.

Applying summation by parts $k$ times starting from (\ref{eq:L-n})
(using every time (\ref{key-ident})), we arrive at the
identity:
\begin{equation}\label{Ln-parts-k}
L_n^{\a,\b}(x) =  c^\diamond
\sum_{j=0}^\infty
A_k(j)\frac{\Gamma(j+\a+k+\b+1)}{\Gamma(j+\b+1)}P_j^{(\a+k,\b)}(x).
\end{equation}
By (\ref{A1}) it readily follows that
\begin{equation}\label{est-A1}
\|A_1^{(m)}\|_{\infty} \le n^{-m-1}\|\wh a^{(m+1)}\|_\infty,
\quad m\ge 0,
\end{equation}
and inductively, using (\ref{Ak+1}), it follows that
\begin{equation}\label{est-Ak}
\|A_k^{(m)}\|_{\infty}
\le cn^{-m-2k+1}\max_{0\le \nu\le m+k}\|\wh a^{(\nu)}\|_\infty,
\quad m\ge 0, \: k\ge 2,
\end{equation}
where $c=c(k,m,\a,\b)$.

Now, from (\ref{est.Pn}) and (\ref{est-Ak}) with $m=0$, we obtain
for $\pi/n \le \theta \le \pi/2$
\begin{align*}
|L_n^{\a,\b}(\cos \theta)|
&\le  c \sum_{j=n-k}^{2n}n^{-2k+1} j^{\a+k-1/2}\theta^{-\a-k-1/2}\\
&\le cn^{2\a+2}(n\theta)^{-\a-k-1/2}
\le c\frac{n^{2\a+2}}{(1+n\theta)^k},
\end{align*}
and for $\pi/2 \le \theta \le \pi$
\begin{align*}
|L_n^{\a,\b}(\cos \theta)|
\le  c \sum_{j=n-k}^{2n}n^{-2k+1} j^{\a+k}j^{\b}
\le cn^{-k+2\a+2}
\le c\frac{n^{2\a+2}}{(1+n\theta)^k}.
\end{align*}
Thus (\ref{est.Ln}) is established when $r=0$.

The case when $r\ge 1$ is an easy consequence of Markov's inequality:
If $Q \in \Pi_m$, then
$
\|Q'\|_{L^\infty[a,b]} \le 2m^2(b-a)^{-1}\|Q\|_{L^\infty[a,b]}.
$

We give the proof for $r=1$ only; in general it follows inductively.
Clearly, (\ref{est.Ln}) with $r=0$ is equivalent to
\begin{equation}\label{est-Lnx}
|L_n^{\a,\b}(x)| \le c\frac{n^{2\a+2}}{(1+n\sqrt{1-x})^k},
\quad x\in [-1,1].
\end{equation}
If $x \in [0,1]$, then by (\ref{est-Lnx}) and Markov's inequality
\begin{align*}
\Big|\frac{d}{dx}L_n^{\a,\b}(x)\Big|
&\le \Big\|\frac{d}{dx}L_n^{\a,\b}\Big\|_{L^\infty[-1,x]}
\le 8n^2(1+x)^{-1} \|L_n^{\a,\b}\|_{L^\infty[-1,x]}\\
&\le c\frac{n^{2\a+4}}{(1+n\sqrt{1-x})^k},
\end{align*}
which is (\ref{est.Ln}) with $r=1$.
For $x \in [-1,0)$ we apply Markov's inequality
on $[-1,0]$ which leads readily to the same result.
The proof of Theorem~\ref{thm:BD} is complete.
\qed


\begin{thebibliography}{99}
\bibitem{Bor}
        J. Borwein,
        A note on Farkas lemma, Utilitas Math. \textbf{24} (1983), 235-241.

\bibitem{BLW}
        L. Bos, N. Levenberg and S. Waldron,
        Metrics associated to multivariate polynomial inequalities,
        in \textit{Advances in Constructive Approximation: Vanderbilt 2003},
        p. 133--147, Nashboro Press, Brentwood, TN, 2004.

\bibitem{BD}
        G. Brown and F. Dai,
        Approximation of smooth functions on compact two-point homogeneous
        spaces,
        J. Funct. Anal. \textbf{220}  (2005),  no. 2, 401--423.


\bibitem{DX}
        C. F. Dunkl and Yuan Xu,
        \textit{Orthogonal polynomials of several variables},
        Cambridge Univ. Press, 2001.

\bibitem{Dz}
        J. Dziuba\'{n}ski,
        Triebel-Lizorkin spaces associated with Laguerre and Hernite expansions,
        Proc. Amer. Math. Soc. \textbf{125} (1997), 3547--3554.

\bibitem{F-J-W}
        M. Frazier, B. Jawerth, and G. Weiss,
        Littlewood-Paley theory and the study of function spaces,
        CBMS No 79 (1991), AMS.


\bibitem{H}
        L. H\"{o}rmander, The spectral function of an elliptic operator,
        Acta Math. \textbf{121} (1968), 193--218.

\bibitem{KPX}
        G. Kyriazis, P. Petrushev, Yuan Xu,
        Jacobi decomposition of weighted Triebel-Lizorkin and Besov spaces, preprint.
        (http://www.math.sc.edu/$\sim$pencho/)

\bibitem{Meyer}
        Y. Meyers,
        {\em Ondelletes et Op\'{e}rateursI: Ondelletes},
        Hermann, Paris, 1990.

\bibitem{MNW1}
        H. Mhaskar, F. Narcowich, and J. Ward,
        Spherical Marcinkiewicz-Zygmund inequalities and positive quadrature,
        \textit{Math. Comp.} \textbf{70} (2001), 1113-1130.

\bibitem{MNW2}
        H. Mhaskar, F. Narcowich, and J. Ward,
        Corrigendum to ``Spherical Marcinkiewicz-Zygmund inequalities
        and positive quadrature''
        \textit{Math. Comp.} \textbf{71} (\#237), pp.~453--454, 2001.

\bibitem{Myso}
        I. P. Mysovskikh
        \textit{Interpolatory cubature formulas},
        in Russian, ``Nauka'', Moscow, 1981.

\bibitem{NPW}
        F. Narcowich, P. Petrushev, and J. Ward,
        Localized tight frames on spheres,
        \textit{SIAM J. Math. Anal.} \textbf{38} (2006), 347-692.

\bibitem{NPW2}
        F. Narcowich, J. Ward, and P. Petrushev,
        Decomposition of Besov and Triebel-Lizorkin spaces on the sphere,
        \textit{J. Funct. Anal.} \textbf{238} (2006), 530-564.

\bibitem{PX}
        P. Petrushev and Yuan Xu,
        Localized polynomial frames on the interval with Jacobi weights,
        \textit{J. Fourier Anal. and Appl.} \textbf{5} (2005), 557--575.

\bibitem{St}
        A. Stroud,
        {\it Approximation calculation of multiple integrals},
        Prentice-Hall, NJ, 1971.

\bibitem{Sz}
        G. Szeg\H{o},
        \textit{Orthogonal Polynomials},
    Amer. Math. Soc. Colloq. Publ. Vol.23, Providence, 4th edition,
        1975.


\bibitem{X98}
        Yuan Xu,
        Orthogonal polynomials and cubature formulae on spheres and on balls,
        \textit{SIAM J. Math. Anal.} \textbf{29} (1998), 779--793.

\bibitem{X99}
        Yuan Xu,
        Summability of Fourier orthogonal series for Jacobi weight
        on a ball in $\RR^d$,
        \textit{Trans. Amer. Math. Soc.} \textbf{351} (1999), 2439-2458.

\bibitem{X04}
        Yuan Xu,
        Weighted approximation of functions on the unit sphere,
        \textit{Const. Approx.} \textbf{21} (2005), 1-28.

\end{thebibliography}
\end{document}